\documentclass{amsart}

\usepackage[english]{babel}
\usepackage[utf8]{inputenc}

\usepackage{amsmath}
\usepackage{amssymb}
\usepackage{amsthm}
\usepackage{mathtools}
\usepackage{mathrsfs}
\usepackage{abraces}
\usepackage{bbm}

\usepackage{tikz}
\usetikzlibrary{cd,patterns,positioning,decorations.pathreplacing}
\usepackage{xstring}
\usepackage{ifthen}
\usepackage{pdflscape}
\usepackage{diagbox}

\usepackage{titlesec}
\titleformat{\section}{\normalfont\fillast\scshape}{\thesection.}{.5em}{}
\titleformat{\subsection}[runin]{\normalfont\bfseries}{\thesubsection.}{.5em}{}[.]

\usepackage{enumitem}

\setlist[description]{%
	font={\normalfont\itshape},
	}
\renewcommand\labelenumi{(\textit{\roman{enumi}})}
\renewcommand\theenumi\labelenumi

\theoremstyle{plain}
\newtheorem{thm}{Theorem}[subsection]
\newtheorem{lem}[thm]{Lemma}
\newtheorem{pro}[thm]{Proposition}
\newtheorem{cor}[thm]{Corollary}
\theoremstyle{definition}
\newtheorem{dfn}[thm]{Definition}
\theoremstyle{remark}
\newtheorem{rmk}[thm]{Remark}
\newtheorem{exa}[thm]{Example}

\newcommand{\Kb}{\mathcal{K}^{\mathrm{b}}}

\newcommand{\uom}{\underline{\omega}}
\newcommand{\uw}{\underline{w}}
\newcommand{\ux}{\underline{x}}

\newcommand{\DbOZ}{\mathcal{D}^{\mathrm{b}}(\mathcal{O}_0^{\mathbb{Z}})}

\DeclareMathOperator{\Hom}{Hom}

\DeclareMathOperator{\id}{id}

\tikzset{%
	pics/dot/.style n args={1}{%
		code={%
			\draw (0,0) -- (0,#1);
			\fill (0,#1) circle (1.5pt); 		
			}
	}
}

\tikzset{%
	positive/.style = {line width=2pt},
	negative/.style = {line width=.25pt,double distance=1.5pt},
    posneg/.style = {line width=2pt, dashed}
}

\tikzset{%
	pics/patch/.style n args={1}{
		code={
			\node[scale=#1,circle,draw,fill=white,inner sep=0,font=\tiny] at (0,0) {$\times$};
			}
	}
}

\title{\textbf{A Spectral Sequence for Rouquier Complexes}}
\author{Leonardo Maltoni}

\begin{document}

\maketitle
\begin{abstract}
    We define a spectral sequence converging to the morphism spaces between Rouquier complexes. From the first page we deduce the cohomological range of morphisms between standard Rouquier complexes. We also extend some results by Ko and Mazorchuk on extension groups between Verma modules to the modular version of (graded) category $\mathcal{O}$ and refine a numerical criterion to detect ``additional'' extensions in terms of $R$-polynomials.
\end{abstract}
\section*{Introduction}
Many categories of interest in representation theory have braid group actions. A classical example,  defined in various ways in \cite{KhomMaz,MazStr,Rou_cat}, is the action on the bounded derived category 
of the principal block of category $\mathcal{O}$. 
The endofunctors induced by this action can be modeled by Rouquier complexes, which provide a categorified version of the braid group \cite{Rou_cat}. These complexes also control other classical braid group actions, see for example \cite{KhoSei}, and play a role in the theory of categorified link invariants, see \cite{Khov}.

Rouquier complexes are defined in terms of the Hecke category $\mathscr{H}$.
This admits various algebraic or geometric descriptions: for example, Soergel bimodules \cite{Abe, Soe_HC, Soe}, sheaves over the flag variety \cite{JMW,KazLus_intcoh,RW,Soe}, moment graphs \cite{BraMacMomentGraphs, Fiebig2008, Fiebig2010}. The works of Elias, Khovanov and Williamson \cite{Elias_two,EliKho,EW} provided a diagrammatic presentation by generators and relations. This version was adapted by the author \cite{Mal_pat} to a diagrammatic description of a dg enhancement of the category of Rouquier complexes.
The main result is that, after introducing a suitable pre-order on the bases for the dg morphism spaces deduced from Libedinsky \textit{light leaves} \cite{Lib}, the differential exhibits an upper-triangularity property.

In this paper, we obtain a spectral sequence $(E_r^{p,q})$ from the filtration induced by the above pre-order. This converges to the actual morphism spaces between Rouquier complexes in the homotopy category $\Kb(\mathscr{H})$ and the first page is easy to compute. From this we obtain the cohomological range of morphisms between \textit{standard} Rouquier complexes, that is, those corresponding to positive lifts of elements in the Coxeter group $W$. Under certain assumptions, if $x,y\in W$ and $F_x,F_y$ are the corresponding Rouquier complexes, then:
\[
    \Hom^i(F_x,F_y)=0 \quad \text{unless $x\le y$ and $d(x,y)\le i \le \ell(y)-\ell(x),$}
\]
where $\le$ denotes the Bruhat order, whereas $\ell(x),\ell(y)$ are the lengths of $x,y$ in $W$ and $d(x,y)$ is the distance between $x$ and $y$ in the Bruhat graph of $W$. More generally one has to consider the dimension of the span of the roots corresponding to the reflections involved in the paths from $x$ to $y$.

Finally we apply our spectral sequence to the graded category $\mathcal{O}$. Over the complex numbers, the classical construction by Beilinson, Ginzburg and Soergel \cite{BeiGinSoe} of the graded derived category $\DbOZ$ and its standard heart $\mathcal{O}_0^{\mathbb{Z}}$ agrees with the homotopy category of Soergel modules and the heart of a suitable perverse t-structure \cite{EWHodge,Soergel1990} (making the above action of Rouquier complexes more explicit). This can be generalised to the homotopy category $\Kb(\overline{\mathscr{H}})$ of the regular module $\overline{\mathscr{H}}$ over the Hecke category, thus providing a modular version for $\DbOZ$: see \cite{AchRic_modII,Makisumi2023}. The morphisms between the images of standard Rouquier complexes then correspond to the extensions between graded Verma modules:
\begin{equation}
\label{eq:extensionRouquier}
    \mathrm{Ext}^k_{\mathcal{O}_0^{\mathbb{Z}}}(\Delta_x\langle i\rangle,\Delta_y)\cong \Hom_{\Kb(\overline{\mathscr{H}})}(\overline{F}_{w_0x^{-1}}\langle i\rangle,\overline{F}_{w_0y^{-1}}[k])
\end{equation}
where on the right hand side $\langle 1 \rangle = [-1](1)$ is the \textit{Tate twist}.

The appropriate version $(\overline{E}_r^{p,q})$ of our spectral sequence converges at the second page to these extensions. More precisely, let $x\ge y$ (otherwise the extension group \eqref{eq:extensionRouquier} is zero) and let $h\in\{0,1\}$ be such that $\ell(x)-\ell(y) \equiv h \mod 2$. The term $\overline{E}^{p,q}_1$ is homogeneous of polynomial degree $h-2q$ (i.e., it is isomorphic to a sum of copies of $\Bbbk(2q-h)$). This implies that the dimension of \eqref{eq:extensionRouquier} is equal to $\dim \overline{E}^{p,q}_2$ when $k=p-q+h$ and $i=h-2q$ (and zero anywhere else). 
However, the differentials $\overline{E}_1^{p,q}\to \overline{E}_1^{p+1,q}$ of the first page involve some complicated diagrammatic calculations and are hard to determine. Nevertheless the first page alone recovers some results by Ko and Mazorchuk \cite{MazKo} and provides a criterion to detect \textit{additional} extensions. 
Consider the $R$-polynomial $R_{x,y}$ and recall there is a polynomial $\widetilde{R}_{x,y}$ with nonnegative integer coefficients such that $R_{x,y}(v)=\widetilde{R}_{x,y}(v-v^{-1})$. Then the Poincaré polynomial $P(s,t):=\sum_{p,q}(\dim \overline{E}_1^{p,q})s^pt^q$ of the first page has the following description:
\begin{equation}
    \label{eq:firstpage}
    P(s,t)=
    s^{-\frac{h}{2}} t^{\frac{h}{2}} \widetilde{R}_{x,y}(s^{\frac{1}{2}}(t^{-\frac{1}{2}}+t^{\frac{1}{2}}))
\end{equation}
If we arrange the potential contributions to \eqref{eq:extensionRouquier} according to $k$ and $i$ we recover the same region of existence as in \cite{MazKo}. The following depicts an example with $\ell(x)-\ell(y)=7$.
\[
\begin{tikzpicture}[yscale=0.8,xscale=0.8,baseline=-1.5cm]
    \foreach \x in {-3,-2,-1,0,1,2,3,4}{
    \node[font=\tiny] at (\x,0) {$E_1^{3,\x}$};
    }
    \foreach \x in {-2,-1,0,1,2,3}{
    \node[font=\tiny]  at (\x,-1) {$E_1^{2,\x}$};
    \draw[-latex] (\x,-0.7)--(\x,-0.3);
    }
    \foreach \x in {-1,0,1,2}{
    \node[font=\tiny] at (\x,-2) {$E_1^{1,\x}$};
    \draw[-latex] (\x,-1.7)--(\x,-1.3);
    }
    \foreach \x in {0,1}{
     \node[font=\tiny] at (\x,-3) {$E_1^{0,\x}$};
    \draw[-latex] (\x,-2.7)--(\x,-2.3);
    }
\end{tikzpicture}
\quad \rightsquigarrow \qquad
\begin{tikzpicture}[scale=.4,baseline=0]
    \draw[red] (0,0)--(7,7)--(0,-7)--cycle;
    \fill[pink] (0,0)--(7,7)--(0,-7)--cycle;
    \draw[gray,->] (-1,0)--(8,0) node[right] {$k$};
    \draw[gray,->] (0,-8)--(0,8) node[above] {$i$};
    \foreach \y in {1,2,3,4,5,6,7,-1,-2,-3,-4,-5,-6,-7}{
    \draw[lightgray] (-.5,\y)--(7.5,\y);
    }
    \foreach \x in {1,2,3,4,5,6,7}{
        \draw[lightgray] (\x,-7.5)--(\x,7.5);
    }
    \fill (0,-1) circle (4pt); \fill[white] (0,-1)circle (3pt); 
    \fill (0,-3) circle (4pt); \fill[white] (0,-3)circle (3pt);
    \fill (0,-5) circle (4pt); \fill[white] (0,-5)circle (3pt);
    \fill (0,-7) circle (4pt);
    \fill (1,1) circle (4pt); \fill[white] (1,1)circle (3pt);
    \fill (1,-1) circle (4pt);
    \fill (1,-3) circle (4pt);
    \fill (1,-5) circle (4pt);
    \fill (2,1) circle (4pt);
    \fill (2,-1) circle (4pt);
    \fill (2,-3) circle (4pt);
    \fill (3,3) circle (4pt); \fill[white] (3,3)circle (3pt);
    \fill (3,1) circle (4pt);
    \fill (3,-1) circle (4pt);
    \fill (4,3) circle (4pt);
    \fill (4,1) circle (4pt);
    \fill (5,5) circle (4pt);\fill[white] (5,5)circle (3pt);
    \fill (5,3) circle (4pt);
    \fill (6,5) circle (4pt);
    \fill (7,7) circle (4pt);
\end{tikzpicture}
\]
Furthermore, consider the increasing alternating sums $S_{r}^q=\sum_{p=0}^r (-1)^p\dim \overline{E}_1^{p,q}$ along each column of the first page. These sequences are obtained from $P(-1,t)$ so they can be easily computed from the $R$-polynomials. If the sequence $(S_0^q,S_1^q,\dots )$ does not alternate in sign, then there is an additional extension.
\section{The Hecke category and Rouquier complexes}
\subsection{Notation}
Let $S$ be a finite set and, for $s,t\in S$, let $m_{st}\in \mathbb{N}\cup\{\infty\}$ such that $m_{ss}=1$ for all $s\in S$. Let $(W,S)$ be the associated Coxeter system:
\begin{equation}
	W=\langle s\in S\mid (st)^{m_{st}}=1\rangle
\end{equation}
Elements of $S$ are then involutions and are called \emph{simple reflections}. Their conjugates are called \emph{reflections}. 

Let $S^*$ be the free monoid generated by $S$ whose elements $\uw$ are expressions in the alphabet $S$. Applying the obvious map $S^*\rightarrow W$, every expression represents an element $w$ of $W$ (in general different expressions can represent the same element). An expression $\uw=s_1s_2\dots s_k$ representing $w$ is called \emph{reduced} if there is no shorter expression for $w$ and in this case we set $\ell(w):=k$.

We set $w\le wt$ when $t$ is a reflection and $\ell(wt)>\ell(w)$. The transitive closure of this relation on $W$, still denoted by $\le$, is the \emph{Bruhat order}.

The Hecke algebra $\mathbf{H}$ associated with $(W,S)$ is the $\mathbb{Z}[v,v^{-1}]$-algebra  generated by elements $\delta_s$, for all $s\in S$, and relations:
\begin{align}
    &\delta_s^2=1-(v-v^{-1})\delta_s,\, \forall s\in S&
    &\underbrace{\delta_s\delta_t\delta_s\dots }_{m_{st}\text{ terms}}= \underbrace{\delta_t\delta_s\delta_t\dots }_{m_{st}\text{ terms}},\, \forall s,t\in S,  \label{eq_braidrel}
\end{align}
where the second relation is omitted when $m_{st}=\infty$.

To the Coxeter system $(W,S)$ one associates the (generalised) braid group, or Artin-Tits group:
\begin{equation}
    B_W:=\langle s \in S \mid \underbrace{sts\dots}_{m_{st}\text{ terms}}=\underbrace{tst\dots}_{m_{st}\text{ terms}}\rangle
\end{equation}
Notice that to simplify notation we see $S$ as generating both $W$ and $B_W$. In particular, when seen in $B_W$, the elements $s$ and $s^{-1}$ are different. Let $\Sigma=S\cup S^{-1}$ and denote by $\Sigma^*$ the free monoid generated by $\Sigma$. An element $\uom\in \Sigma^*$ is called \emph{braid expression} and it represents an element $\omega\in B_W$ via the natural map $\Sigma^*\rightarrow B_W$.

\subsection{Some formulas in the Hecke algebra}
For any expression $\underline{w}=s_1\cdots s_k\in S^*$ let $\delta_{\underline{w}}:=\delta_{s_1}\cdots \delta_{s_k}$. Now let $x\in W$ and let $\underline{x}$ be a reduced expression for $x$. One can prove that the definition $\delta_x:=\delta_{\underline{x}}$ is independent of the choice of a reduced expression. Recall that the elements $\delta_x$ for $x\in W$ form a basis of $\mathbf{H}$ over $\mathbb{Z}[v,v^{-1}]$ called the \textit{standard basis}. 

For every $s\in S$ the element $\delta_s$ is invertible in $\mathbf{H}$ with $\delta_s^{-1}=\delta_s + (v-v^{-1})$. Hence one can define an involution  $\overline{(\cdot)}: \mathbf{H}\to \mathbf{H}$ via the assignments $v\mapsto v^{-1}$ and $\delta_s\mapsto\delta_s^{-1}$. So for $\uw=s_1\dots s_k\in S^*$ we have $\overline{\delta_{\uw}}=\delta_{s_1}^{-1}\cdots \delta_{s_k}^{-1}$ and for $x\in W$ we have $\overline{\delta_x}=(\delta_{x^{-1}})^{-1}$, which defines a \textit{dual} basis.

We want to give an expression for $\overline{\delta_{\underline{w}}}$ in terms of the standard basis.
One can easily prove the following formula:
\begin{align}
    \label{eq: product formulas in the Hecke algebra}
    & \delta_x\delta_s^{-1}=\begin{cases}
        \delta_{xs} + (v-v^{-1})\delta_x & \text{if $xs>x$,} \\
        \delta_{xs} & \text{if $xs<x$.}
    \end{cases}
\end{align}
From this one can then recursively compute the decompositions of arbitrary products $\overline{\delta_{\underline{w}}}$ into the basis $\{\delta_x\}_{x\in W}$, but there is also a closed formula.
		\begin{dfn}
			Given $\uw=s_1\dots s_k \in S^*$, an element $\eta\in\{0,1\}^k$ is called a \emph{01-sequence} for $\uw$. Associated to a 01-sequence $\epsilon$ one defines a \emph{Bruhat stroll} $x_0$, $x_1$, \dots, $x_k$, of elements of $W$, where:
			\begin{align}			
				&x_0:=1,&
				&x_i:=s_1^{\eta_1}\dots s_i^{\eta_i},& \text{for } i=1,\dots, k.
			\end{align}		
			where, for $s\in S$, we set $s^0=1$.
			One says that $\eta$ \emph{expresses} the element $x$, and one writes $\underline{w}^\eta = x$, if $x_k=x$.	 
		\end{dfn}
		We can encode some information about the Bruhat stroll induced by a 01-sequence $\eta$ for $\uw$ by \textit{decorating} it as follows. Each entry $\eta_i$ is decorated by:
		\begin{itemize}
			\item a $U$, if $x_{i}s_{i+1}>x_i$ in the Bruhat order;
			\item a $D$, if $x_{i}s_{i+1}<x_{i}$.
		\end{itemize}
\begin{dfn}
    A 01-sequence $\eta\in\{0,1\}^k$ for $\underline{w}=s_1\dots s_k$ is called \textit{distinguished} if its decoration has no $D0$. The set of disringuished 01-sequences for $\uw$ is denoted by $\mathcal{D}(\uw)$.
\end{dfn}
For every distinguished 01-sequence $\eta$ let $g(\eta)$ be the number of 0's (hence $U0$'s) that it contains.
Then one obtains the following.
\begin{pro}
    \label{pro: closed formula for product of inverses}
    Let $\uw=s_1\dots s_k\in S^*$ be any expression. Then:
    \begin{equation}
        \label{eq: closed formula for product of inverses}
        \overline{\delta_{\uw}} = \sum_{\eta \in \mathcal{D}(\underline{w})} (v-v^{-1})^{g(\eta)}\delta_{\underline{w}^\eta}
    \end{equation}
\end{pro}
The change of basis between $\{\delta_x\}_{x\in W}$ and $\{\overline{\delta_x}\}_{x\in W}$ is given by the $R$-polynomials:
\begin{align}
    \label{eq: change of basis standard - dual}
    &\overline{\delta_{y}}=\sum_{x\in W}R_{x,y}(v^{-1})\delta_x&
    &\delta_{y}=\sum_{x\in W}R_{x,y}(v)\overline{\delta_x}
\end{align}
Hence Proposition \ref{pro: closed formula for product of inverses} implies the following classical properties of $R$-polynomials: see \cite[\S 5.1-\S 5.3]{BjoBre}.
\begin{pro}
\label{pro: classical properties of R - polynomials}
Let $x,y\in W$. We have:
    \begin{enumerate}
        \item The polynomial $R_{x,y}(v)$ is 0 unless $x\le y$.
        \item We have $R_{x,y}(v)=\widetilde{R}_{x,y}(v-v^{-1})$, where $\widetilde{R}_{x,y}(u)=\sum_{\eta} u^{g(\eta)}$, the sum begin over $\eta\in \mathcal{C}(\underline{y})$ with $\underline{y}^{\eta}=x$.
    \end{enumerate}
\end{pro}
\begin{proof}
    Both parts follow from equation \eqref{eq: closed formula for product of inverses}.
    For part (i) it is sufficient to notice that if $\underline{y}$ is a reduced expression for $y$ then, for all 01-sequence $\eta$ for $\underline{y}$, one has $\underline{y}^{\eta}\le y$. Part (ii) is just a combination of \eqref{eq: closed formula for product of inverses} and \eqref{eq: change of basis standard - dual}.
\end{proof}
Equivalently, one can consider the pairing 
$\beta:\mathbf{H}\times \mathbf{H}\to \mathbb{Z}[v,v^{-1}]$ defined by:
\begin{equation}
    \beta(\delta_x,\overline{\delta_y}):=\begin{cases}
        1 & \text{if $x=y$,}\\
        0 & \text{otherwise,}
    \end{cases}
\end{equation}
and extended $\mathbb{Z}[v,v^{-1}]$-semlinearly. This means that for $a,b\in \mathbf{H}$, one has $\beta(a,bv)=v\beta(a,b)$ and $\beta(va,b)=v^{-1}\beta(a,b)$. Then:
\begin{equation}
    \label{eq: R-polynomial as pairing}
     \beta(\delta_x,\delta_y)=R_{x,y}(v).
\end{equation}
\begin{rmk}
    The Hecke algebra also has the \textit{canonical} or \textit{Kazhdan-Lusztig basis} $\{b_x\}_{x\in W}$, which is self-dual $\overline{b_x}=b_x$ and such that
    $b_y=\sum_{x\le y} h_{x,y}(v) \delta_x$ with $h_{y,y}=1$ and $h_{x,y}(v)\in v\mathbb{Z}[v]$ for $x<y$ (and it is uniquely determined by these properties). The $h_{x,y}$ are (a renormalisation of) the \textit{Kazhdan-Lusztig polynomials}. See \cite{KazLus_can,Soergel1997}.
\end{rmk}
\subsection{Realisations of Coxeter systems}
Let $\Bbbk$ be a commutative ring. The following can be seen as a generalisation of both the geometric representation of Coxeter systems and the natural representations of Weyl groups arising from root data. See \cite[\S 3.1]{EW}.
\begin{dfn}
	A \emph{realisation} of $(W,S)$ over $\Bbbk$ is a free finite rank $\Bbbk$-module $\mathfrak{h}$ equipped with subsets of \emph{simple coroots} and \emph{simple roots}:
	\begin{align}
	&\{\alpha_s^\vee\}_{s\in S}\subset \mathfrak{h},&
	&\{\alpha_s\}_{s\in S}\subset\mathfrak{h}^*:=\Hom(\mathfrak{h},\Bbbk),
	\end{align}
	such that:
	\begin{enumerate}
	\item for all $s\in S$, we have $\langle\alpha_s,\alpha_s^\vee\rangle = 2$;
	\item the action $s(v)=v-\langle \alpha_s,v\rangle\alpha_s^\vee$, for all $s\in S$ and $v\in \mathfrak{h}$, defines a $\Bbbk$-linear representation of $W$ on $\mathfrak{h}$ (and hence also a contragredient representation given by $s(\lambda)=\lambda-\langle\lambda,\alpha_s^\vee\rangle\alpha_s$, for $s\in S$ and $\lambda\in \mathfrak{h}^*$);
	\item the condition \cite[\S 3.1 (3.3)]{EW} holds
	\end{enumerate}
\end{dfn}
\begin{rmk}
The latter is a technical condition which ensures that certain structural morphisms in the Hecke category are well-defined in this generality. It is automatic for root data and geometric representations.
\end{rmk}
Let $R=S(\mathfrak{h}^*)$ be the symmetric algebra over $\mathfrak{h}^*$ with $\deg\mathfrak{h}^*=2$. The linear action of $W$ on $\mathfrak{h}^*$ induces an action of $W$ on $R$ by algebra automorphisms. For each $s\in S$, let $\partial_s:R\to R$ be the \emph{Demazure operator} defined by:
\begin{align}
    &\partial_s(f):=\dfrac{f-s(f)}{\alpha_s}
\end{align}
\begin{dfn}
A realisation $\mathfrak{h}$ of $(W,S)$ over $\Bbbk$ is \emph{admissible} if:
\begin{enumerate}
	\item it satisfies \emph{Demazure surjectivity}: for each $s\in S$ there exists $\delta_s\in \mathfrak{h}^*$ such that $\langle\delta_s,\alpha_s^\vee\rangle=1$;
	\item it is \emph{balanced} in the sense of \cite[Def.\ 3.7]{EW}.
	\item the ring $\Bbbk$ is an integral domain.
\end{enumerate}
\end{dfn}
From now on we will assume that all realisations are admissible. 
\begin{rmk}
	These very mild assumptions are satisfied by the geometric representation and by all the realisations arising from root data over any integral domain (with the caveat that in a few cases it is necessary to assume that 2 is invertible). 
    Yet, they are sufficient to make the diagrammatic Hecke category, introduced in the next section, well-behaved.
	A stronger assumption is needed to make the classical category of Soergel bimodules to work properly. Suppose $\Bbbk$ is a field. A realisation $\mathfrak{h}$ over $\Bbbk$ is \emph{reflection faithful} if the representation of $W$ on $\mathfrak{h}$ is faithful and it induces a bijection between reflections in $W$ and hyperplanes of $\mathfrak{h}$ which are fixed by some element of $W$: see \cite{Soe}. This assumption is satisfied for the geometric representation but fails in many realisations arising from root data for affine Weyl groups.
\end{rmk}
\subsection{The Hecke category}
        \label{sec: Hecke category}
		We will now review the definition of the diagrammatic version of the Hecke category.
		
		Let $\mathfrak{h}$ be an admissible realisation over $\Bbbk$ of the Coxeter system $(W,S)$. We will proceed in two steps. The first is defining the Bott-Samelson category. We will omit several details in the definition, and refer the reader \cite{EW}.
		\begin{dfn}
			The \emph{Bott-Samelson category} $\mathscr{H}_{\mathrm{BS}}$ is the $\Bbbk$-linear monoidal category, enriched in graded $R$-bimodules, defined by generators and relations as follows.
			\begin{enumerate}
				\item The objects are formal tensor products of generators $B_s$ for $s\in S$. So a general object has the form $B_{\uw}:=B_{s_1}\otimes \cdots \otimes B_{s_k}$ for an expression $\uw=s_1\dots s_k$. Let $\mathbf{1}$ denote the monoidal unit. 
				\item We associate a color to each simple reflection.
					The space $\Hom_{\mathscr{H}_{\mathrm{BS}}}(B_{\underline{x}},B_{\uw})$ consists of $\Bbbk$-linear combinations of isotopy classes of \emph{Soergel diagrams}, which are defined as follows.
					A Soergel diagram is a colored, \emph{decorated} planar graph contained in the strip $\mathbb{R}\times [0,1]$, with boundary in $\mathbb{R}\times\{0,1\}$ such that:
					\begin{itemize}
						\item The bottom boundary in $\mathbb{R}\times\{0\}$ (respectively, the top boundary in $\mathbb{R}\times\{1\}$) is the sequence of boundary points colored according to the expression $\underline{x}$ (respectively $\uw$).
						\item The edges of the graph are colored in such a way that those connected with the boundary have the same colors as the boundary points.
						\item The other vertices of the graph are either: 
							\begin{enumerate}
								\item univalent (called \emph{dots}), or;
								\item trivalent	with three edges of the same color, or;
								\item $2m_{st}$-valent with edges of alternating colors corresponding to $s$ and $t$, when $m_{st}<\infty$. They are also called $(s,t)$-ars.
							\end{enumerate}
							\begin{center}
								\begin{tikzpicture}[baseline=0,scale=.7]
									\def\shift{4cm}
									\draw[gray,dashed] (0,0) circle (1cm); %
										\draw[red] (0,-1)--(0,0); \fill[red] (0,0) circle (2pt);
									\node[anchor=north,inner sep=8pt] at (0,-1) {(a) \textit{Dot}};
									\begin{scope}[xshift=\shift]
										\draw[gray,dashed] (0,0) circle (1cm); %
											\draw[red] (0,1)--(0,0); %
												\draw[red] (-30:1cm)--(0,0)--(-150:1cm);
										\node[anchor=north,inner sep=8pt] at (0,-1) {(b) \textit{Trivalent}};
									\end{scope}
									\begin{scope}[xshift=2*\shift,rotate=-15]
										\draw[gray,dashed] (0,0) circle (1cm); %
											\draw[red] (0,1)--(0,0); %
												\draw[red] (-30:1cm)--(150:1cm);%
													\draw[red](-150:1cm)--(30:1cm);
											\draw[blue,rotate=30] (0,1)--(0,0); %
												\draw[blue,rotate=30] (-30:1cm)--(150:1cm);%
													\draw[blue,rotate=30](-150:1cm)--(30:1cm);
										\foreach \i in {-60,-75,-90}{%
												\fill (\i:.5cm) circle (.5pt);
											}
									\end{scope}
									\node[anchor=north,xshift=1.4*\shift,inner sep=8pt] at (0,-1) {(c) \textit{$(s,t)$-ar}};
								\end{tikzpicture}
							\end{center}
						\item Decorations are boxes labeled by homogeneous elements in $R$ that can appear in any region (i.e., connected component of the complement of the graph).
						\item The degree of a Soergel graph is determined by summing the contributions of all its parts: dots contribute 1, trivalent vertices contribute $-1$ and stars contribute 0. The contribution of polynomial boxes is the degree of the corresponding element of $R$.
					\end{itemize}
					Then, composition of morphisms is given by gluing diagrams vertically, whereas tensor product is given by gluing them horizontally. The identity morphism of an object $B_{\uw}$ is the diagram with parallel vertical strands colored according to the expression $\uw$.
				\item These diagrams are identified via the following relations (apart from isotopy):
					\begin{description}
						\item[Polynomial relations.] (involving polynomial boxes)
							\begin{align}
								& \text{\emph{Addition}:}  &
								\begin{tikzpicture}[baseline=-0.1cm,scale=.7, transform shape]
									\draw[gray,dashed] (0,0) circle (1cm); \node[draw] at (0,0) {$f$};
								\end{tikzpicture}
								&+
								\begin{tikzpicture}[baseline=-0.1cm,scale=.7, transform shape]
									\draw[gray,dashed] (0,0) circle (1cm); \node[draw] at (0,0) {$g$};
								\end{tikzpicture}
								=
								\begin{tikzpicture}[baseline=-0.1cm,scale=.7, transform shape]
									\draw[gray,dashed] (0,0) circle (1cm); \node[draw] at (0,0) {$f+g$};
								\end{tikzpicture}
								\\
								& \text{\emph{Multiplication}:}  &
								\begin{tikzpicture}[baseline=-0.1cm,scale=.7, transform shape]
									\draw[gray,dashed] (0,0) circle (1cm); \node[draw] at (-0.4,0.4) {$f$};\node[draw] at (0.4,-0.4) {$g$};
								\end{tikzpicture}
								&=
								\begin{tikzpicture}[baseline=-0.1cm,scale=.7, transform shape]
									\draw[gray,dashed] (0,0) circle (1cm); \node[draw] at (0,0) {$fg$};
								\end{tikzpicture}
                            \end{align}
                            \begin{align}
								& \text{\emph{Barbell} relation:}  &
								\begin{tikzpicture}[baseline=-0.1cm,scale=.7, transform shape] \label{barbell}
									\draw[gray,dashed] (0,0) circle (1cm); \draw[red] (0,.5)--(0,-.5); \fill[red] (0,.5) circle (2pt);\fill[red] (0,-.5) circle (2pt);
								\end{tikzpicture}
								&=
								\begin{tikzpicture}[baseline=-0.1cm,scale=.7, transform shape]
									\draw[gray,dashed] (0,0) circle (1cm); \node[draw] at (0,0) {$\alpha_{\textcolor{red}{s}}$};
								\end{tikzpicture}			\\
								& \text{\emph{Sliding} relation:} & 
								\begin{tikzpicture}[baseline=-0.1cm,scale=.7, transform shape]\label{sliding}
									\draw[gray,dashed] (0,0) circle (1cm); \draw[red] (0,1)--(0,-1); 
									\node[draw] at (.5,0) {$f$};
								\end{tikzpicture}			
								&=
								\begin{tikzpicture}[baseline=-0.1cm,scale=.7, transform shape]
									\draw[gray,dashed] (0,0) circle (1cm); \draw[red] (0,1)--(0,-1); 
									\node[draw, inner sep =.08cm] at (-.5,0) {$\textcolor{red}{s}(f)$};
								\end{tikzpicture}			
								+
								\begin{tikzpicture}[baseline=-0.1cm,scale=.7, transform shape]
									\draw[gray,dashed] (0,0) circle (1cm); \draw[red] (0,1)--(0,.5); \fill[red] (0,.5) circle (2pt);
									\draw[red] (0,-1)--(0,-.5); \fill[red] (0,-.5) circle (2pt); 
									\node[draw] at (0,0) {$\partial_{\textcolor{red}{s}}(f)$};
								\end{tikzpicture}
								\end{align}
						\item[One color relations.] (involving only strands of a single color)
							\begin{align}
								& \text{\emph{Frobenius associativity}:} &
								\begin{tikzpicture}[baseline=-0.1cm,color=red,scale=sqrt(2)/3,scale=.7, transform shape]\label{associativity}
									\draw (-1.5,-1.5)--(-0.6,0)--(0.6,0)--(1.5,-1.5);
									\draw (-1.5,1.5)--(-0.6,0); \draw(0.6,0)--(1.5,1.5);
									\draw[dashed,gray] (0,0) circle ({3/sqrt(2)});
								\end{tikzpicture}
								&=
								\begin{tikzpicture}[baseline=-0.1cm,color=red,scale=sqrt(2)/3,rotate=90,scale=.7, transform shape]
									\draw (-1.5,-1.5)--(-0.6,0)--(0.6,0)--(1.5,-1.5);
									\draw (-1.5,1.5)--(-0.6,0); \draw(0.6,0)--(1.5,1.5);
									\draw[dashed,gray] (0,0) circle ({3/sqrt(2)});
								\end{tikzpicture} \\
								& \text{\emph{Frobenius unit}:} &
								\begin{tikzpicture}[color=red,baseline=-0.1cm,scale=.7, transform shape]\label{dotline}
									\draw[gray,dashed] (0,0) circle (1cm);
									\draw (0,-1) -- (0,1); \draw (0,0)--(.3,0);\fill (0.3,0) circle (2pt);
								\end{tikzpicture}
								&=
								\begin{tikzpicture}[baseline=-0.1cm,scale=.7, transform shape]
									\draw[gray,dashed] (0,0) circle (1cm);						
									\draw[red] (0,-1) -- (0,1);
								\end{tikzpicture}
								\\
								&\text{\emph{Needle} relation:} &
								\begin{tikzpicture}[baseline=0,scale=.7, transform shape]\label{eq_needle}
									\draw[gray,dashed] (0,0) circle (1cm);
									\draw[red] (0,-1)--(0,-0.3) arc (-90:270:0.3cm);
								\end{tikzpicture}
								&= 0
							\end{align}
						\item[Two color relations.] (involving stars: two versions are given according to the parity of $m_{st}$)
							\begin{equation}\label{2mtriv}
								\begin{tikzpicture}[baseline=0,scale=.8,transform shape]
									\draw[dashed,gray] (0,0) circle (1cm);
									\coordinate (a) at (-.3,0); \coordinate (b) at (0.5,0);
									\draw[red] (180:1cm)--(a); \draw[blue] ({180-360/11}:1cm)--(a);
									\draw[red] ({180-2*360/11}:1cm)--(a); \draw[blue] ({180-5*360/11}:1cm)--(b);
									\draw[red] ({180-4*360/11}:1cm)--(a);
									\draw[blue] ({180+360/11}:1cm)--(a);
									\draw[red] ({180+2*360/11}:1cm)--(a); \draw[blue] ({180+5*360/11}:1cm)--(b);
									\draw[red] ({180+4*360/11}:1cm)--(a);
									\draw[blue](a)--(b);
									\node at (-0.05,.5) {...};\node at (-0.05,-.5) {...};
								\end{tikzpicture}
								=
								\begin{tikzpicture}[baseline=0,scale=.8,transform shape]
									\draw[dashed,gray] (0,0) circle (1cm);
									\coordinate (a) at (-.5,0); \coordinate (b) at (0.15,0.35); \coordinate (c) at (0.15,-0.35);
									\draw[red] (180:1cm)--(a); \draw[blue] ({180-360/11}:1cm)--(b);
									\draw[red] ({180-2*360/11}:1cm)--(b); \draw[blue] ({180-5*360/11}:1cm)--(b);
									\draw[red] ({180-4*360/11}:1cm)--(b);
									\draw[blue] ({180+360/11}:1cm)--(c);
									\draw[red] ({180+2*360/11}:1cm)--(c); \draw[blue] ({180+5*360/11}:1cm)--(c);
									\draw[red] ({180+4*360/11}:1cm)--(c);
									\draw[blue](c)..controls (-0.15,0)..(b); \draw[red] (b) -- (c);\draw[red] (b)--(a)--(c);
									\draw[red] (b) ..controls (0.7,0).. (c);
									\node at (.2,.6) {...};\node at (.2,-.6) {...};\node at (.35,0) {...};
								\end{tikzpicture}\quad\text{or}\quad
								\begin{tikzpicture}[baseline=0,scale=.8,transform shape]
									\draw[dashed,gray] (0,0) circle (1cm);
									\coordinate (a) at (-.3,0); \coordinate (b) at (0.5,0);
									\draw[red] (180:1cm)--(a); \draw[blue] ({180-360/11}:1cm)--(a);
									\draw[red] ({180-2*360/11}:1cm)--(a); \draw[red] ({180-5*360/11}:1cm)--(b);
									\draw[blue] ({180-4*360/11}:1cm)--(a);
									\draw[blue] ({180+360/11}:1cm)--(a);
									\draw[red] ({180+2*360/11}:1cm)--(a); \draw[red] ({180+5*360/11}:1cm)--(b);
									\draw[blue] ({180+4*360/11}:1cm)--(a);
									\draw[red](a)--(b);
									\node at (-0.05,.5) {...};\node at (-0.05,-.5) {...};
								\end{tikzpicture}
								=
								\begin{tikzpicture}[baseline=0,scale=.8,transform shape]
									\draw[dashed,gray] (0,0) circle (1cm);
									\coordinate (a) at (-.5,0); \coordinate (b) at (0.15,0.35); \coordinate (c) at (0.15,-0.35);
									\draw[red] (180:1cm)--(a); \draw[blue] ({180-360/11}:1cm)--(b);
									\draw[red] ({180-2*360/11}:1cm)--(b); \draw[red] ({180-5*360/11}:1cm)--(b);
									\draw[blue] ({180-4*360/11}:1cm)--(b);
									\draw[blue] ({180+360/11}:1cm)--(c);
									\draw[red] ({180+2*360/11}:1cm)--(c); \draw[red] ({180+5*360/11}:1cm)--(c);
									\draw[blue] ({180+4*360/11}:1cm)--(c);
									\draw[blue](c)..controls (-0.15,0)..(b); \draw[red] (b) -- (c);\draw[red] (b)--(a)--(c);
									\draw[blue] (b) ..controls (0.7,0).. (c);
									\node at (.2,.6) {...};\node at (.2,-.6) {...};\node at (.35,0) {...};
								\end{tikzpicture}
							\end{equation}
							and
							\begin{equation}\label{eq_2mdot}
								\begin{tikzpicture}[baseline=0,scale=.8,transform shape]
									\draw[dashed,gray] (0,0) circle (1cm);
									\coordinate (a) at (-.2,0); \coordinate (b) at (0.5,0);
									\draw[red] (180:1cm)--(a); \draw[blue] ({180-360/11}:1cm)--(a);
									\draw[blue] ({180-3*360/11}:1cm)--(a); 
									\draw[red] ({180-4*360/11}:1cm)--(a);
									\draw[blue] ({180+360/11}:1cm)--(a);
									\draw[blue] ({180+3*360/11}:1cm)--(a); \fill[blue] (b) circle (2pt);
									\draw[red] ({180+4*360/11}:1cm)--(a);
									\draw[blue](a)--(b);
									\node at (-0.4,.5) {...};\node at (-0.4,-.5) {...};
								\end{tikzpicture}
								=
								\begin{tikzpicture}[baseline=0,scale=.8,transform shape]
									\draw[dashed,gray] (0,0) circle (1cm);
									\coordinate (a) at (-.2,0); \coordinate (b) at (0.5,0);
									\draw[red] (180:1cm)--(a); \draw[blue] ({180-360/11}:1cm)--(a);
									\draw[blue] ({180-3*360/11}:1cm)--(a); 
									\draw[blue] ({180+360/11}:1cm)--(a);
									\draw[blue] ({180+3*360/11}:1cm)--(a);
									\draw[red] ({180+4*360/11}:1cm) -- (b) -- ({180-4*360/11}:1cm); \draw[red](a)--(b);
									\node at (-0.4,.6) {...};\node at (-0.4,-.6) {...};
									\node[draw,circle, inner sep=.1cm,fill=white] at (a) {$\mathrm{JW}$};
								\end{tikzpicture}
								\quad\text{or}\quad
								\begin{tikzpicture}[baseline=0,scale=.8,transform shape]
									\draw[dashed,gray] (0,0) circle (1cm);
									\coordinate (a) at (-.2,0); \coordinate (b) at (0.5,0);
									\draw[red] (180:1cm)--(a); \draw[blue] ({180-360/11}:1cm)--(a);
									\draw[red] ({180-3*360/11}:1cm)--(a); 
									\draw[blue] ({180-4*360/11}:1cm)--(a);
									\draw[blue] ({180+360/11}:1cm)--(a);
									\draw[red] ({180+3*360/11}:1cm)--(a); \fill[red] (b) circle (2pt);
									\draw[blue] ({180+4*360/11}:1cm)--(a);
									\draw[red](a)--(b);
									\node at (-0.4,.5) {...};\node at (-0.4,-.5) {...};
								\end{tikzpicture}
								=
								\begin{tikzpicture}[baseline=0,scale=.8,transform shape]
									\draw[dashed,gray] (0,0) circle (1cm);
									\coordinate (a) at (-.2,0); \coordinate (b) at (0.5,0);
									\draw[red] (180:1cm)--(a); \draw[blue] ({180-360/11}:1cm)--(a);
									\draw[red] ({180-3*360/11}:1cm)--(a); 
									\draw[blue] ({180+360/11}:1cm)--(a);
									\draw[red] ({180+3*360/11}:1cm)--(a);
									\draw[blue] ({180+4*360/11}:1cm) -- (b) -- ({180-4*360/11}:1cm); \draw[blue](a)--(b);
									\node at (-0.4,.6) {...};\node at (-0.4,-.6) {...};
									\node[draw,circle, inner sep=.1cm,fill=white] at (a) {$\mathrm{JW}$};
								\end{tikzpicture}
							\end{equation}
							Here the $\mathrm{JW}$-circles are the \emph{Jones-Wenzel morphisms} defined in \cite[\S\,5.2]{EW}.
						\item[Three color relations.] We will omit these relations and point the the reader to \cite[\S\ 1.4.3]{EW}.
					\end{description}				
			\end{enumerate}
   		This concludes the definition of $\mathscr{H}_{\mathrm{BS}}$.
		\end{dfn}
         Then one defines the full Hecke category in a completely formal way from $\mathscr{H}_{\mathrm{BS}}$.
		\begin{dfn}
			The \emph{diagrammatic Hecke category} $\mathscr{H}$ is the Karoubi envelope of the closure of $\mathscr{H}_{\mathrm{BS}}$
			by shifts and direct sums.
		\end{dfn}
        The split Grothendieck group $K_{\oplus}(\mathscr{H})$ of the Hecke category is naturally a $\mathbb{Z}[v,v^{-1}]$-algebra, with $[B_1][B_2]:=[B_1\otimes B_2]$ and $v[B]:=[B(1)]$.
        \begin{thm}
            \label{thm: categorification}
            Let $\Bbbk$ be a complete local ring. Then the map $K_{\oplus}(\mathscr{H})\to\mathbf{H}$ uniquely determined by $[B_s]\mapsto \delta_s+v$, is an isomorphism. Furthermore the morphism spaces $\Hom_{\mathscr{H}}(B_1,B_2)$ are free as graded left (or right) $R$-modules and the graded dimension is given by:
            \begin{equation}
                \label{eq: graded dimension of morphism space}
                \mathrm{grdim}_R\Hom(B_1,B_2) = \beta([B_1],[B_2])
            \end{equation}
        \end{thm}
        In the Theorem, the graded dimension of a free $R$-module is encoded in a Laurent polynomial in the following way. If:
        \begin{equation}
            \label{eq: definition of graded dimension}
            M\cong\bigoplus_{i\in \mathbb{Z}} R(-i)^{\oplus a_i},
        \end{equation}
        then we set $\mathrm{grdim}_R( M):=\sum_i a_i v^i\in \mathbb{Z}[v,v^{-1}]$. In other words the coefficient $a_i$ is the number of elements of degree $i$ in a basis for $M$.
		\begin{rmk}
			Under the same assumptions of the Theorem, one can classify the indecomposable objects of $\mathscr{H}$. For every $x\in W$ there is an indecomposable object $B_x$, which can be characterised as the top summand of any Bott-Samelson $B_{\underline{x}}$ for any reduced expression $\underline{x}$ for $x$. It satisfies $[B_x]=b_x$ and, up to shift, all indecomposable objects are of this form.
		\end{rmk}
		\begin{rmk}
			When the realisation $\mathfrak{h}$ is reflection faithful this category is equivalent to the classical category of Soergel bimodules. There is another version of the latter category defined by Abe \cite{Abe} which restores the equivalence also when the reflection faithfulness assumption is dropped.
			There are also other versions of the Hecke category, equivalent to the present one under appropriate assumptions, as mentioned in the introduction.
		\end{rmk}
\subsection{Light leaves maps}
		The morphism spaces in the Hecke category are free finite rank graded left (or right) $R$-modules. One can describe explicit bases for the morphism spaces between Bott-Samelson objects via the \textit{light leaves} of Libedinsky \cite{Lib}. For our purposes we will only need to deal with the case of morphisms towards the unit object:
		\begin{equation}\label{eq_homsp}
			\Hom_{\mathscr{H}}(B_{\uw},\mathbf{1}).
		\end{equation}
        A basis for \eqref{eq_homsp} is labeled by 01-sequences for $\uw$ expressing the identity element.
        For every such 01-sequence $\eta$, the corresponding basis morphism $L_\eta: B_{\uw}\rightarrow \mathbf{1}$ is constructed as follows, according to some choices.
		\begin{dfn}
            \label{def: recursive procedure for light leaves}
			First, for each $i<k$, choose a reduced expression $\ux_i$ for $x_i$ and hence a corresponding object $B_i$. In particular $B_0=\mathbf{1}$. We define $L_\eta$ by induction. Let $\lambda_0$ to be $\id_{\mathbf{1}}$. Then, suppose we have constructed a morphism $\lambda_i:B \rightarrow B_i$, and consider the decorated value of $\eta_{i+1}$.
			\begin{enumerate}
				\item If $\eta_{i+1}=U0$, then $x_i=x_{i+1}$ and we choose a \emph{braid move} $\Phi$ (i.e., a sequence of braid relations) from $\ux_i$ to $\ux_{i+1}$.
				\item If $\eta_{i+1}=U1$, then $x_{i+1}=x_is_{i+1}$ and we choose a braid move $\Phi$ from the reduced expression $\ux_i s_{i+1}$ to $\ux_{i+1}$.
				\item If $\eta_{i+1}=D0$, then there is a reduced expression for $x_{i}$ of the form $\underline{z} s_{i+1}$ where $\underline{z}$ is a reduced expression for $z=x_{i}s_{i+1}$, and we choose a braid move $\Psi$ from $\ux_i$ to $\underline{z} s_{i+1}$. Furthermore, we choose another braid move $\Phi$ from the latter to $\ux_{i+1}$.
				\item Finally, if $\eta_{i+1}=D1$, then we choose, as in the previous case, a braid move $\Psi$ from $\ux_{i}$ to $\underline{z} s_{i+1}$. This time $\underline{z}$ is a reduced expression for $x_{i+1}$ and we choose another braid move $\Phi$ from $\underline{z}$ to $\ux_{i+1}$.
			\end{enumerate}
			Then $\lambda_{i+1}$ is defined by one of the following diagrams, according to the decorated value of $\eta_{i+1}$ (here we denote $\beta(\Phi)$ and $\beta(\Psi)$ the compositions of $(s,t)$-ars corresponding to the braid moves $\Phi$ and $\Psi$):
			\[
				\begin{tikzpicture}[scale=.7]
					\draw[gray] (-.5,0)--(3,0); 
					\draw[violet] (0,0)--(0,.5);\draw[violet] (0.5,0)--(0.5,.5);\draw[violet] (1.5,0)--(1.5,.5);\node at (1,.25) {\dots};\draw[violet] (2,0)--(2,.5);
					\draw (-.25,.5)--(2.25,.5)--(2.25,1.5)--(-.25,1.5)--cycle; \node at (1,1) {$\lambda_i$};
					\draw[violet] (0,1.5)--(0,2);\draw[violet] (0.5,1.5)--(0.5,2);\draw[violet] (1.5,1.5)--(1.5,2);\node at (1,1.75) {\dots};\draw[violet] (2,1.5)--(2,2);
					\draw[red] (2.5,0)--(2.5,1);\fill[red](2.5,1)circle (2pt);
					\draw (-.25,2)--(2.25,2)--(2.25,2.5)--(-.25,2.5)--cycle; \node at (1,2.25) {\footnotesize $\beta(\Phi)$};
					\draw[violet] (0,2.5)--(0,3);\draw[violet] (0.5,2.5)--(0.5,3);\draw[violet] (1.5,2.5)--(1.5,3);\node at (1,2.75) {\dots};\draw[violet] (2,2.5)--(2,3);
					\draw[gray] (-.5,3)--(3,3);
					\node[font=\small] at (1.25,-.5){$\eta_{i+1}=U0$};
				\end{tikzpicture}\qquad
				\begin{tikzpicture}[scale=.7]
					\draw[gray] (-.5,0)--(3,0); 
					\draw[violet] (0,0)--(0,.5);\draw[violet] (0.5,0)--(0.5,.5);\draw[violet] (1.5,0)--(1.5,.5);\node at (1,.25) {\dots};\draw[violet] (2,0)--(2,.5);
					\draw (-.25,.5)--(2.25,.5)--(2.25,1.5)--(-.25,1.5)--cycle; \node at (1,1) {$\lambda_i$};
					\draw[violet] (0,1.5)--(0,2);\draw[violet] (0.5,1.5)--(0.5,2);\draw[violet] (1.5,1.5)--(1.5,2);\node at (1,1.75) {\dots};\draw[violet] (2,1.5)--(2,2);
					\draw[red] (2.5,0)--(2.5,2);
					\draw (-.25,2)--(2.75,2)--(2.75,2.5)--(-.25,2.5)--cycle; \node at (1.25,2.25) {\footnotesize $\beta(\Phi)$};
					\draw[violet] (0,2.5)--(0,3);\draw[violet] (0.5,2.5)--(0.5,3);\draw[violet] (1.5,2.5)--(1.5,3);\node at (1,2.75) {\dots};\draw[violet] (2,2.5)--(2,3);\draw[violet] (2.5,2.5)--(2.5,3);
					\draw[gray] (-.5,3)--(3,3);
					\node[font=\small] at (1.25,-.5){$\eta_{i+1}=U1$};
				\end{tikzpicture}\qquad
				\begin{tikzpicture}[scale=.7]
					\draw[gray] (-.5,0)--(3,0); 
					\draw[violet] (0,0)--(0,.5);\draw[violet] (0.5,0)--(0.5,.5);\draw[violet] (1.5,0)--(1.5,.5);\node at (1,.25) {\dots};\draw[violet] (2,0)--(2,.5);
					\draw (-.25,.5)--(2.25,.5)--(2.25,.9)--(-.25,.9)--cycle; \node at (1,.7) {\footnotesize $\lambda_i$};
					\draw[violet] (0,.9)--(0,1.1);\draw[violet] (0.5,0.9)--(0.5,1.1);\draw[violet] (1.5,0.9)--(1.5,1.1);\draw[violet] (2,0.9)--(2,1.1);\node at (1,1) {\dots};
					\draw (-.25,1.1)--(2.25,1.1)--(2.25,1.5)--(-.25,1.5)--cycle; \node at (1,1.3) {\footnotesize $\beta(\Psi)$};
					\draw[violet] (0,1.5)--(0,2);\draw[violet] (0.5,1.5)--(0.5,2);\draw[red] (2,1.5) arc (180:0:.25cm);\node at (1,1.75) {\dots};\draw[violet] (1.5,1.5)--(1.5,2);
					\draw[red] (2.5,0)--(2.5,1.5); \draw[red] (2.25,2)--(2.25,1.75);
					\draw (-.25,2)--(2.5,2)--(2.5,2.5)--(-.25,2.5)--cycle; \node at (1.125,2.25) {\footnotesize $\beta(\Phi)$};
					\draw[violet] (0,2.5)--(0,3);\draw[violet] (0.5,2.5)--(0.5,3);\node at (1.125,2.75) {\dots};\draw[violet] (1.75,2.5)--(1.75,3);\draw[violet] (2.25,2.5)--(2.25,3);
					\draw[gray] (-.5,3)--(3,3);
					\node[font=\small] at (1.25,-.5){$\eta_{i+1}=D0$};
				\end{tikzpicture}\qquad
				\begin{tikzpicture}[scale=.7]
					\draw[gray] (-.5,0)--(3,0); 
					\draw[violet] (0,0)--(0,.5);\draw[violet] (0.5,0)--(0.5,.5);\draw[violet] (1.5,0)--(1.5,.5);\node at (1,.25) {\dots};\draw[violet] (2,0)--(2,.5);
					\draw (-.25,.5)--(2.25,.5)--(2.25,.9)--(-.25,.9)--cycle; \node at (1,.7) {\footnotesize $\lambda_i$};
					\draw[violet] (0,.9)--(0,1.1);\draw[violet] (0.5,0.9)--(0.5,1.1);\draw[violet] (1.5,0.9)--(1.5,1.1);\draw[violet] (2,0.9)--(2,1.1);\node at (1,1) {\dots};
					\draw (-.25,1.1)--(2.25,1.1)--(2.25,1.5)--(-.25,1.5)--cycle; \node at (1,1.3) {\footnotesize $\beta(\Psi)$};
					\draw[violet] (0,1.5)--(0,2);\draw[violet] (0.5,1.5)--(0.5,2);\draw[red] (2,1.5) arc (180:0:.25cm);\node at (1,1.75) {\dots};\draw[violet] (1.5,1.5)--(1.5,2);
					\draw[red] (2.5,0)--(2.5,1.5);
					\draw (-.25,2)--(1.75,2)--(1.75,2.5)--(-.25,2.5)--cycle; \node at (0.75,2.25) {\footnotesize $\beta(\Phi)$};
					\draw[violet] (0,2.5)--(0,3);\draw[violet] (0.5,2.5)--(0.5,3);\node at (1,2.75) {\dots};\draw[violet] (1.5,2.5)--(1.5,3);
					\draw[gray] (-.5,3)--(3,3);
					\node[font=\small] at (1.25,-.5){$\eta_{i+1}=D1$};
				\end{tikzpicture}
				\]
				Then one can check that the target of $\lambda_{i+1}$ corresponds with the reduced expression $\ux_{i+1}$ and the induction can continue.
				Finally we define $L_\eta$ to be $\lambda_k \in \Hom_{\mathscr{H}}(B,\mathbf{1})$.
		\end{dfn}			
		\begin{thm}[see \cite{EW,Lib}]
			For a given Bott-Samelson object $B_{\uw}$ in $\mathscr{H}_{\mathrm{BS}}$, for any set of choices as above, the set $\{L_{\eta}\mid \uw^{\eta}=1 \}$ is a basis for $\Hom_{\mathscr{H}}(B_{\uw},\mathbf{1})$ as a left $R$-module.
		\end{thm}
\subsection{Rouquier complexes}
		Consider, for each simple reflection $s\in S$, the following \emph{standard} and \emph{costandard} complexes:
		\[
			\begin{tikzcd}[row sep=small,%
							execute at end picture={%
								\def\hdott{.32cm}
								\pic[red,yshift=.3em] at (A) {dot={\hdott}};
								\pic[red,yshift=1.5em,yscale=-1] at (B) {dot={\hdott}};
								}]
				F_{\textcolor{red}{s}}=\,\cdots\ar[r]		& 0 \ar[r]		& 0	\ar[r]								& \underline{B_{\textcolor{red}{s}}}\ar[r,""{coordinate,name=A}]	& \mathbf{1}(1)\ar[r]	& 0\ar[r]		& \cdots \\
				F_{\textcolor{red}{s}}^{-1}=\,\cdots\ar[r]	& 0\ar[r]		& \mathbf{1}(-1)\ar[r,""{coordinate,name=B}]	& \underline{B_{\textcolor{red}{s}}}\ar[r]							& 0	\ar[r]		& 0\ar[r]		& \cdots 
			\end{tikzcd}
		\]
		where, here and in the sequel, the underlined object is understood to be in cohomological degree 0.
		\begin{dfn}		
			If $\uom=s_1^{e_1}\dots s_k^{e_k}\in \Sigma^*$ is a braid expression (with $e_i\in\{\pm 1\}$), we put:
			\begin{equation}\label{eq_stdrou}
				F_{\uom}^\bullet:=F_{s_1}^{e_1}\otimes \dots \otimes F_{s_k}^{e_k}.
			\end{equation}
			Objects of this form are called	\emph{Rouquier complexes}.
		\end{dfn}
		\begin{pro}\label{pro_braidRou}
			In the homotopy category $\mathcal{K}(\mathscr{H})$, we have:
			\begin{enumerate}
			 	\item Let $s\in S$, then $F_s \otimes F_s^{-1}\cong F_s^{-1}\otimes F_s\cong \mathbf{1}[0]$.
				\item Let $s,t\in S$ with $m_{st}<\infty$, then:
					\[
						\underbrace{F_s\otimes F_t\otimes F_s \otimes \cdots}_{m_{st} \text{ terms}} \cong %
							\underbrace{F_t \otimes F_s \otimes F_t\otimes \cdots}_{m_{st} \text{ terms}}
					\]
				\item (Rouquier Canonicity)
					For each pair $\uom_1,\uom_2\in \Sigma^*$ expressing the same element in $B_W$, we have:
					\[
						\Hom_{\mathcal{K}(\mathscr{H})}(F_{\uom_1}^\bullet,F_{\uom_2}^\bullet)\cong R,
					\]
					and one can chose homotopy equivalences $\gamma_{\uom_1}^{\uom_2}\in \Hom(F_{\uom_1}^\bullet,F_{\uom_2}^\bullet)$ such that the system $\{\gamma_{\uom_1}^{\uom_2}\}_{\uom_1,\uom_2}$ is transitive, namely:
					\begin{align*}
						\gamma_{\uom}^{\uom}=&\id_{F^\bullet_{\uom}},& 										&\text{for all $\uom$,}\\
						\gamma_{\uom_2}^{\uom_3}\circ\gamma_{\uom_1}^{\uom_2}=&\gamma_{\uom_1}^{\uom_3},&	&\text{for all $\uom_1$, $\uom_2$, $\uom_3$.}						
					\end{align*}
			\end{enumerate}
		\end{pro}
		Taking limits of these transitive systems, one defines , for $\omega\in B_W$, objects $F_\omega$ in $\mathcal{K}(\mathscr{H})$, unique up to canonical isomorphism, so that $F_\omega\otimes F_{\omega'}=F_{\omega\omega'}$. These give a strict version of Rouquier complexes.

        Recall that the triangulated Grothendieck group of $\Kb(\mathscr{H})$ naturally coincides with the split Grothendieck group of $\mathscr{H}$:
        \begin{equation}
            \label{eq: equality of Grothendieck groups}
            K_{\Delta}(\Kb(\mathscr{H}))=K_{\oplus}(\mathscr{H}),
        \end{equation}
        where $[F^\bullet]=\sum_{i\in \mathbb{Z}} (-1)^{i}[F_i]$. Hence we have $[F_s]=\delta_s$ and $[F_s^{-1}]=\delta_s^{-1}$ so that:
        \begin{equation}
            \label{eq: decategorification of Rouquier complexes}
            [F_{\uom}^\bullet] = \delta_{s_1}^{e_1}\cdots \delta_{s_k}^{e_k}
        \end{equation}
\section{A spectral sequence}
\subsection{Review of the dg setting}
The expansion of the tensor product \eqref{eq_stdrou}
is easy to describe. Each term of the expansion is determined by picking either $\mathbf{1}$ or $B_{s_i}$ in the various factors. We encode this as follows.
\begin{dfn}
Given $\uom=s_1^{e_1}\dots s_k^{e_k}\in \Sigma^*$, a subset $\alpha \subset \{1,\dots, k\}$ is called \emph{arrangement} for $\uom$.
\end{dfn}
Now fix a braid expression $\uom=s_1^{e_1}\dots s_k^{e_k}$. The elements of an arrangement will represent the indices for which we picked the unit. For any arrangement $\alpha$, let:
\begin{align}
	&B_\alpha:= A_1\otimes A_2\otimes \dots \otimes A_k, &
	&\text{where}&
	A_i:=\begin{cases}
            B_{s_i} & \text{if $i\notin \alpha$,} \\
            \mathbf{1} & \text{if $i\in \alpha$,}
    \end{cases}	
\end{align}
Then the following is immediate from the definition of tensor product of complexes.
\begin{pro}
As a graded object, the complex \eqref{eq_stdrou} is:
\begin{align}\label{eq_decomp}
    &\bigoplus_{\alpha\subset \{1,\dots k\}} B_{\alpha}[-q_\alpha](q_\alpha)
    &\text{with}&
    &q_\alpha=\sum_{i\in \alpha}e_i
\end{align}
The components $d_\alpha^{\alpha'}:B_\alpha(q_\alpha)\rightarrow B_{\alpha'}(q_{\alpha'})$ of the differential are non-zero only if $\alpha'$ is obtained from $\alpha$ by adding or removing exactly one index $i$. In this case:
\begin{align}\label{eq_diffcomp}
 	&d_\alpha^{\alpha'}=(-1)^{p(i)}	 	
	\begin{tikzpicture}[x=0.3cm,y=0.4cm,baseline=0.36cm,every node/.style={font=\scriptsize}]
		\draw[gray] (0.5,0)--(9.5,0);\draw[gray] (0.5,2)--(9.5,2);
		\foreach \i in {1,2,4,6,8,9}{%
				\draw[violet] (\i,0)--(\i,2);
			}
		\draw[violet] (5,0)--(5,.6); \fill[violet] (5,.6) circle (1.5pt);
		\node[font=\normalsize] at (3,1) {\dots}; \node[font=\normalsize] at (7,1) {\dots};
	\end{tikzpicture}&
	&\text{or}&
	&d_\alpha^{\alpha'}=(-1)^{p(i)}
		\begin{tikzpicture}[x=0.3cm,y=0.4cm,baseline=0.36cm,every node/.style={font=\scriptsize}]
			\draw[gray] (0.5,0)--(9.5,0);\draw[gray] (0.5,2)--(9.5,2);
			\foreach \i in {1,2,4,6,8,9}{%
					\draw[violet] (\i,0)--(\i,2);
				}
			\draw[violet] (5,2)--(5,1.4); \fill[violet] (5,1.4) circle (1.5pt);
			\node[font=\normalsize] at (3,1) {\dots}; \node[font=\normalsize] at (7,1) {\dots};
		\end{tikzpicture},
\end{align}
where $p(i)$ is the number of $j<i$ with $j\in\alpha$.
\end{pro} 
We consider the above complexes in the dg category $\mathcal{C}_{\mathrm{dg}}(\mathscr{H})$ where the morphism spaces are themselves complexes:
\begin{equation}
	\Hom^\bullet(F^\bullet,G^\bullet):=\bigoplus_{i\in \mathbb{Z}}\Hom^i(F^\bullet, G^\bullet):= \bigoplus_{i\in \mathbb{Z}}\prod_{j\in \mathbb{Z}} \Hom_{\mathscr{H}}(F^j,G^{j+i}),
\end{equation}
with differential $d:\Hom^i(F^\bullet, G^\bullet)\rightarrow \Hom^{i+1}(F^\bullet, G^\bullet)$ defined by:
\begin{equation}\label{eq_diffdg}
	d((f^j)_{j\in \mathbb{Z}}) = (d_{G}^{j+i}\circ f^{j} - (-1)^i f^{j+1}\circ d_F^{j})_{j\in \mathbb{Z}}.
\end{equation}
An easy check gives the following well-known formula.
\begin{lem}
The cohomology of the dg morphism complex recovers the morphism spaces in the homotopy category:
\begin{equation}
	H^i(\Hom^\bullet(F^\bullet, G^\bullet))=\Hom_{\mathcal{K}(\mathscr{H})}(F^\bullet, G^\bullet [i])
\end{equation}
\end{lem}
Now we will recall the formalism of \cite{Mal_pat} to describe these dg morphism complexes in terms of Soergel diagrams with patches. We will restrict to the case
\begin{equation}\label{eq_morphto1}
    \Hom^\bullet(F_{\uom},\mathbf{1}).
\end{equation}
One can indeed reduce to this case via adjunction.
\subsection{Morphisms towards the unit} Thanks to \eqref{eq_decomp}, the morphism space \eqref{eq_morphto1} decomposes as:
\begin{equation}\label{eq_decompmor}
    \bigoplus_{\alpha\subset \{1,\dots,k\}}\Hom(B_\alpha,\mathbf{1})[q_\alpha](-q_\alpha),
\end{equation}
where each summand is now a morphism space in the Hecke category and can be described by Soergel diagrams. This allows us to adapt the definition of Soergel diagram to the dg setting.
\begin{dfn}
	A \textit{Soergel diagram with patches} with source $\uom$ (and empty target) and arrangement $\alpha$, is obtained as follows:
	\begin{itemize}
		\item consider the planar strip $\mathbb{R}\times [0,1]$;
		\item mark on the bottom boundary line the sequence of colors corresponding to the letters of $\uom$.
		\item to keep track of the sign, the segment of the boundary line containing the marking for $s_i$ will be black if $e_i=+1$ and white if $e_i=-1$;
		\item cover with a ``patch'' the markings for which $i\in \alpha$;
		\item draw a standard Soergel diagram only involving the markings not covered by a patch.
	\end{itemize}	 
\end{dfn}
\begin{pro}
Elements of \eqref{eq_morphto1} are linear combinations of Soergel diagrams with patches with source $\uom$, and those with arrangement $\alpha$ belong to the corresponding summand in \eqref{eq_decompmor}, in particular ther cohomological degree is $q_\alpha$.
\end{pro}
By combining \eqref{eq_diffcomp} and \eqref{eq_diffdg}, one obtains a diagrammatic description of the differential in terms of the two operations of \emph{sprouting} and \emph{uprooting}:
\begin{align*}
	&S_i:\quad\begin{tikzpicture}[baseline=.2cm]
		\draw[positive] (0,0)--(1,0);
		\pic[violet] at (.5,0) {patch={1}};
		\node[below=.1cm] at (.5,0) {$i$};
	\end{tikzpicture}
	\rightsquigarrow
	\begin{tikzpicture}[baseline=.2cm]
        \draw[violet] (.5,0)--(.5,.8);
        \fill[violet] (.5,.8) circle (2pt);
		\draw[positive] (0,0)--(1,0);
		\node[below=.1cm] at (.5,0) {$i$};
	\end{tikzpicture}&
	&U_i:\quad \begin{tikzpicture}[baseline=.2cm]
		\draw[violet] (.5,0)--(.5,.5);
		\draw[violet,dashed] (.5,.5)--(.5,1);
		\draw[negative] (0,0)--(1,0);
		\node[below=.1cm] at (.5,0) {$i$};
	\end{tikzpicture}
	\rightsquigarrow
	\begin{tikzpicture}[baseline=.2cm]
		\fill[violet] (.5,.3) circle (2pt);
		\draw[violet] (.5,.3)--(.5,.5);
		\draw[violet,dashed] (.5,.5)--(.5,1);
		\draw[negative] (0,0)--(1,0);
		\pic[violet] at (.5,0) {patch={1}};
		\node[below=.1cm] at (.5,0) {$i$};
	\end{tikzpicture}
\end{align*}					 	
These are to be considered local operations on a diagram with patches.
\begin{pro}
The image via the differential map of a diagram with patches $D$, with source $\uom$ and arrangement $\alpha$, is:
\begin{equation}
	d(D)=\sum_{\substack{ i \in \alpha \\ e_i=+1}} (-1)^{p(i)+q_\alpha+1} S_i(D) + \sum_{\substack{i\notin\alpha \\ e_i=-1}}(-1)^{p(i)+q_\alpha+1} U_i(D)
\end{equation}
where $p(i)$ is again the number of patches to the left of $i$.
\end{pro}
\subsection{Bases of morphisms} 
As shown in \cite[\S 3.6]{Mal_pat}, the dg morphism spaces of the form
$\Hom^\bullet(F_{\uom},\mathbf{1)}$ admit bases of \emph{dg light leaves maps}. These are obtained by first choosing a Bott-Samelson factor in the source complex (which amounts to selecting an arrangement of patches on the expression of $\uom$) and then choosing a light leaves map from this to $\mathbf{1}$ in the classical sense. This can be encoded as follows.
\begin{dfn}
    Given a braid expression $\uom=s_1^{e_1}\dots s_k^{e_k}\in \Sigma^*$, and a $01$-sequence $\eta\in \{0,1\}^{k}$, one defined the corresponding Bruhat stroll in the same way as in the Coxeter case, ignoring the exponents $e_i$.
    \begin{align}
        &x_0:=1,&       &x_i:=s_1^{\eta_1}\dots s_i^{\eta_i}, &   &\text{for $i=1,\dots, k$}.
    \end{align}
    We say that $\eta$ \textit{expresses the unit}, and write $\uom^\eta=1$, if $x_k=1$.
\end{dfn}
\begin{dfn}
    Given a braid expression $\uom$ and an arrangement $\alpha\subset\{1,\dots, k\}$, a 01-sequence $\eta\in \{0,1\}^k$ is said to be \textit{compatible with $\alpha$} if for all $i\in \alpha$ one has $\eta_i=0$.
\end{dfn}
A basis for $\Hom^\bullet(F_{\uom},\mathbf{1)}$ is labeled by pairs $(\alpha,\eta)$ where $\alpha\subset\{1,\dots k\}$ is an arrangement and $\eta\in \{0,1\}^k$ is a 01-sequence compatible with $\alpha$.
\begin{dfn}
    The \emph{dg light leaves map} $L_{\alpha,\eta}$ is the dg morphism from $F_{\uom}$ to $\mathbf{1}$ whose Soergel diagram with patches has source $\uom$ and arrangement $\alpha$ and it is obtained with the same exact procedure as in Definition  \ref{def: recursive procedure for light leaves}, except that, when constructing $\lambda_{i+1}$ from $\lambda_i$, the possibilities are:
    \begin{itemize}
        \item If $i+1\notin \alpha$ then we have the same four cases:
			\[
				\begin{tikzpicture}[scale=.7]
					\draw[violet] (0,0)--(0,.5);\draw[violet] (0.5,0)--(0.5,.5);\draw[violet] (1.5,0)--(1.5,.5);\node at (1,.25) {\dots};\draw[violet] (2,0)--(2,.5);
                    \draw[negative] (-.5,0)--(3,0); 
					\draw[positive, dashed] (-.5,0)--(3,0); 
					\draw (-.25,.5)--(2.25,.5)--(2.25,1.5)--(-.25,1.5)--cycle; \node at (1,1) {$\lambda_i$};
					\draw[violet] (0,1.5)--(0,2);\draw[violet] (0.5,1.5)--(0.5,2);\draw[violet] (1.5,1.5)--(1.5,2);\node at (1,1.75) {\dots};\draw[violet] (2,1.5)--(2,2);
					\draw[red] (2.5,0)--(2.5,1);\fill[red](2.5,1)circle (2pt);
					\draw (-.25,2)--(2.25,2)--(2.25,2.5)--(-.25,2.5)--cycle; \node at (1,2.25) {\footnotesize $\beta(\Phi)$};
					\draw[violet] (0,2.5)--(0,3);\draw[violet] (0.5,2.5)--(0.5,3);\draw[violet] (1.5,2.5)--(1.5,3);\node at (1,2.75) {\dots};\draw[violet] (2,2.5)--(2,3);
					\draw (-.5,3)--(3,3);
					\node[font=\small] at (1.25,-.5){$\eta_{i+1}=U0$};
				\end{tikzpicture}\qquad
				\begin{tikzpicture}[scale=.7]
					\draw[violet] (0,0)--(0,.5);\draw[violet] (0.5,0)--(0.5,.5);\draw[violet] (1.5,0)--(1.5,.5);\node at (1,.25) {\dots};\draw[violet] (2,0)--(2,.5);
                    \draw[negative] (-.5,0)--(3,0); 
					\draw[positive, dashed] (-.5,0)--(3,0); 
					\draw (-.25,.5)--(2.25,.5)--(2.25,1.5)--(-.25,1.5)--cycle; \node at (1,1) {$\lambda_i$};
					\draw[violet] (0,1.5)--(0,2);\draw[violet] (0.5,1.5)--(0.5,2);\draw[violet] (1.5,1.5)--(1.5,2);\node at (1,1.75) {\dots};\draw[violet] (2,1.5)--(2,2);
					\draw[red] (2.5,0)--(2.5,2);
					\draw (-.25,2)--(2.75,2)--(2.75,2.5)--(-.25,2.5)--cycle; \node at (1.25,2.25) {\footnotesize $\beta(\Phi)$};
					\draw[violet] (0,2.5)--(0,3);\draw[violet] (0.5,2.5)--(0.5,3);\draw[violet] (1.5,2.5)--(1.5,3);\node at (1,2.75) {\dots};\draw[violet] (2,2.5)--(2,3);\draw[violet] (2.5,2.5)--(2.5,3);
					\draw[gray] (-.5,3)--(3,3);
					\node[font=\small] at (1.25,-.5){$\eta_{i+1}=U1$};
				\end{tikzpicture}\qquad
				\begin{tikzpicture}[scale=.7]
					\draw[violet] (0,0)--(0,.5);\draw[violet] (0.5,0)--(0.5,.5);\draw[violet] (1.5,0)--(1.5,.5);\node at (1,.25) {\dots};\draw[violet] (2,0)--(2,.5);
                    \draw[negative] (-.5,0)--(3,0); 
					\draw[positive, dashed] (-.5,0)--(3,0); 
					\draw (-.25,.5)--(2.25,.5)--(2.25,.9)--(-.25,.9)--cycle; \node at (1,.7) {\footnotesize $\lambda_i$};
					\draw[violet] (0,.9)--(0,1.1);\draw[violet] (0.5,0.9)--(0.5,1.1);\draw[violet] (1.5,0.9)--(1.5,1.1);\draw[violet] (2,0.9)--(2,1.1);\node at (1,1) {\dots};
					\draw (-.25,1.1)--(2.25,1.1)--(2.25,1.5)--(-.25,1.5)--cycle; \node at (1,1.3) {\footnotesize $\beta(\Psi)$};
					\draw[violet] (0,1.5)--(0,2);\draw[violet] (0.5,1.5)--(0.5,2);\draw[red] (2,1.5) arc (180:0:.25cm);\node at (1,1.75) {\dots};\draw[violet] (1.5,1.5)--(1.5,2);
					\draw[red] (2.5,0)--(2.5,1.5); \draw[red] (2.25,2)--(2.25,1.75);
					\draw (-.25,2)--(2.5,2)--(2.5,2.5)--(-.25,2.5)--cycle; \node at (1.125,2.25) {\footnotesize $\beta(\Phi)$};
					\draw[violet] (0,2.5)--(0,3);\draw[violet] (0.5,2.5)--(0.5,3);\node at (1.125,2.75) {\dots};\draw[violet] (1.75,2.5)--(1.75,3);\draw[violet] (2.25,2.5)--(2.25,3);
					\draw[gray] (-.5,3)--(3,3);
					\node[font=\small] at (1.25,-.5){$\eta_{i+1}=D0$};
				\end{tikzpicture}\qquad
				\begin{tikzpicture}[scale=.7]
					\draw[violet] (0,0)--(0,.5);\draw[violet] (0.5,0)--(0.5,.5);\draw[violet] (1.5,0)--(1.5,.5);\node at (1,.25) {\dots};\draw[violet] (2,0)--(2,.5);
                    \draw[negative] (-.5,0)--(3,0); 
					\draw[positive, dashed] (-.5,0)--(3,0); 
					\draw (-.25,.5)--(2.25,.5)--(2.25,.9)--(-.25,.9)--cycle; \node at (1,.7) {\footnotesize $\lambda_i$};
					\draw[violet] (0,.9)--(0,1.1);\draw[violet] (0.5,0.9)--(0.5,1.1);\draw[violet] (1.5,0.9)--(1.5,1.1);\draw[violet] (2,0.9)--(2,1.1);\node at (1,1) {\dots};
					\draw (-.25,1.1)--(2.25,1.1)--(2.25,1.5)--(-.25,1.5)--cycle; \node at (1,1.3) {\footnotesize $\beta(\Psi)$};
					\draw[violet] (0,1.5)--(0,2);\draw[violet] (0.5,1.5)--(0.5,2);\draw[red] (2,1.5) arc (180:0:.25cm);\node at (1,1.75) {\dots};\draw[violet] (1.5,1.5)--(1.5,2);
					\draw[red] (2.5,0)--(2.5,1.5);
					\draw (-.25,2)--(1.75,2)--(1.75,2.5)--(-.25,2.5)--cycle; \node at (0.75,2.25) {\footnotesize $\beta(\Phi)$};
					\draw[violet] (0,2.5)--(0,3);\draw[violet] (0.5,2.5)--(0.5,3);\node at (1,2.75) {\dots};\draw[violet] (1.5,2.5)--(1.5,3);
					\draw[gray] (-.5,3)--(3,3);
					\node[font=\small] at (1.25,-.5){$\eta_{i+1}=D1$};
				\end{tikzpicture}
				\]
                \item If $i+1\in \alpha$ then $\eta_{i+1}$ must be 0 (by the compatibility requirement) and, a part from braid moves, we only cover the $(i+1)$-th boundary point with a patch.   
                \[
				\begin{tikzpicture}[scale=.7]
					\draw[violet] (0,0)--(0,.5);\draw[violet] (0.5,0)--(0.5,.5);\draw[violet] (1.5,0)--(1.5,.5);\node at (1,.25) {\dots};\draw[violet] (2,0)--(2,.5);
                    \draw[negative] (-.5,0)--(3,0); 
					\draw[positive, dashed] (-.5,0)--(3,0); 
					\draw (-.25,.5)--(2.25,.5)--(2.25,1.5)--(-.25,1.5)--cycle; \node at (1,1) {$\lambda_i$};
					\draw[violet] (0,1.5)--(0,2);\draw[violet] (0.5,1.5)--(0.5,2);\draw[violet] (1.5,1.5)--(1.5,2);\node at (1,1.75) {\dots};\draw[violet] (2,1.5)--(2,2);
					\draw (-.25,2)--(2.25,2)--(2.25,2.5)--(-.25,2.5)--cycle; \node at (1,2.25) {\footnotesize $\beta(\Phi)$};
					\draw[violet] (0,2.5)--(0,3);\draw[violet] (0.5,2.5)--(0.5,3);\draw[violet] (1.5,2.5)--(1.5,3);\node at (1,2.75) {\dots};\draw[violet] (2,2.5)--(2,3);
					\draw[gray] (-.5,3)--(3,3);
                    \pic[red] at (2.5,0) {patch={1}};
					\node[font=\small] at (1.25,-.5){$\eta_{i+1}=U0$};
				\end{tikzpicture}\qquad
				\begin{tikzpicture}[scale=.7]
					\draw[violet] (0,0)--(0,.5);\draw[violet] (0.5,0)--(0.5,.5);\draw[violet] (1.5,0)--(1.5,.5);\node at (1,.25) {\dots};\draw[violet] (2,0)--(2,.5);
                    \draw[negative] (-.5,0)--(3,0); 
					\draw[positive, dashed] (-.5,0)--(3,0); 
					\draw (-.25,.5)--(2.25,.5)--(2.25,.9)--(-.25,.9)--cycle; \node at (1,.7) {\footnotesize $\lambda_i$};
					\draw[violet] (0,.9)--(0,1.1);\draw[violet] (0.5,0.9)--(0.5,1.1);\draw[violet] (1.5,0.9)--(1.5,1.1);\draw[violet] (2,0.9)--(2,1.1);\node at (1,1) {\dots};
					\draw (-.25,1.1)--(2.25,1.1)--(2.25,1.5)--(-.25,1.5)--cycle; \node at (1,1.3) {\footnotesize $\beta(\Psi)$};
					\draw[violet] (0,1.5)--(0,2);\draw[violet] (0.5,1.5)--(0.5,2);
                    \node at (1,1.75) {\dots};\draw[violet] (1.5,1.5)--(1.5,2);
					\draw[red] (2,2)--(2,1.5);
					\draw (-.25,2)--(2.25,2)--(2.25,2.5)--(-.25,2.5)--cycle; \node at (1,2.25) {\footnotesize $\beta(\Phi)$};
					\draw[violet] (0,2.5)--(0,3);\draw[violet] (0.5,2.5)--(0.5,3);\node at (1,2.75) {\dots};\draw[violet] (1.5,2.5)--(1.5,3);\draw[violet] (2,2.5)--(2,3);
					\draw[gray] (-.5,3)--(3,3);
					\node[font=\small] at (1.25,-.5){$\eta_{i+1}=D0$};
                    \pic[red] at (2.5,0) {patch={1}};
				\end{tikzpicture}
                \]
        \end{itemize}
\end{dfn}
\begin{thm}
    The set $\{L_{\alpha,\eta}\mid \uom^\eta=1\}$ is a basis for the dg $R$-module $\Hom^\bullet(F_{\uom},\mathbf{1})$.
\end{thm}
\begin{rmk}    
Notice that the Bruhat stroll (and hence the decorations for $\eta$) only depend on the sequence $s_1\dots s_k$ and are independent of the exponents $e_1,\dots, e_k$. These nevertheless are important to determine both the (cohomological and polynomial) degree of $L_{\alpha,\eta}$ and its behaviour under the differential map.
\end{rmk}
\subsection{Upper triangularity of the differential}
One of the main results of \cite{Mal_pat} is a description of the behaviour of the differential with respect to the dg light leaves basis. More precisely it is upper triangular with respect to the \emph{path-dominance (pre-)order}:
\begin{dfn}
    Given 01-sequences $\eta$ and $\theta$ for $\uom$, with Bruhat strolls $x_0,\dots x_k$ and $y_0,\dots,y_k$ respectively, we set $\eta\lesssim\theta$ if, for all $0\le i\le k$, we have $x_i\le y_i$ in the Bruhat order.
\end{dfn}
This defines a partial order on 01-sequence \cite{EW}, and hence a preorder on pairs $(\alpha,\eta)$, where two pairs are equivalent when the second entry is the same.
Then the upper triangularity result is the following (see \cite[Prop.\ 3.10]{Mal_pat}).
\begin{thm}\label{thm_unitr}
    The differential is upper triangular with respect to the path-dominance pre-order. More precisely, if:
    \[
        d(L_{\alpha,\eta}) = \sum_{\gamma, \theta}a_{\gamma,\theta} L_{\gamma,\theta},
    \]
    for $a_{\gamma,\theta}\in R$, then $a_{\gamma,\theta}\neq 0$ only if $\theta\lesssim\eta$.
\end{thm}
\subsection{Path-dominance filtration and spectral sequence}\label{sec: spectral sequence}
The result in the previous section induces a filtration in the complex $\Hom^\bullet(F_{\uom},\mathbf{1})$ which can be used to define a spectral sequence that converges to:
\[
    H^*\big(\Hom^\bullet(F_{\uom},\mathbf{1})\big)=\Hom^*(F_\omega,\mathbf{1}).
\]
We will now construct precisely such a filtration.
\begin{dfn}
    A 01-sequence $\eta$ for $\uom=s_1^{e_1}\cdots s_k^{e_k}\in \Sigma^*$ expressing the unit is called \textit{critical} if it has no negative $D0$ or positive $U0$. In other words, for all $i$, if $\eta_i=U0$ then $e_i=-1$ and if $\eta_i=D0$ then $e_i=-1$. The set of critical 01-sequences for $\uom$ is denoted by $\mathcal{C}(\uom)$
\end{dfn}
\begin{dfn}
    Let $\uom=\sigma_1\dots \sigma_k$ be a braid word. A function $f:\mathcal{C}(\uom)\to \mathbb{Z}_{>0}$ is called \textit{path-dominant} if for any $\eta,\theta\in \mathcal{C}(\uom)$ such that $\eta \lesssim \theta$ we have $f(\eta)\ge f(\theta)$ with the equality holding if and only if $\eta=\theta$.
\end{dfn}
\begin{exa}\label{ex: path-dominant functions}
    Let $W$ be any Coxeter group. Here are two example of path-dominant functions.
    \begin{enumerate}
        \item Fix any $\uom\in \Sigma^*$. For $\eta$ a 01-sequences for $\uom$ with Bruhat strolls $x_1,\dots, x_k$, set $f(\eta)=\ell(x_1)+\dots+\ell(x_k)$. Then it follows immediately from the definition that $f$ is a path-dominant function (the inequality is respected for all pairs of 01-sequences, not just critical ones).
        \item Let $\uw$ be a reduced expression in $S^*$ and let $\uom$ be its negative lift to $\Sigma^*$, that is if $\uw=s_1\dots s_k$ then $\uom=s_1^{-1}\dots s_k^{-1}$. Then notice that critical 01-sequences for $\uom$ are the same as distinguished 01-sequences expressing the unit. Let $h\in \{0,1\}$ such that $h\equiv k\mod 2$. For every such $\eta$ let $g(\eta)$ be the number of 0s and define $f(\eta):=(g(\eta)-h)/2$. This gives a path-dominant function that we call \textit{standard}.
    \end{enumerate}
\end{exa}
Now, fix a path-dominant function $f$ for $\uom$. Let $C_0=\Hom^\bullet(F_{\uom},\mathbf{1})$. For $p> 0$, consider the pairs $(\alpha,\eta)$ with the property that there exist a $\theta\in \mathcal{C}(\uom)$ with $\eta\lesssim \theta$ and $f(\theta)\ge p$. Let $C_p$ be the dg $R$-submodule of $\Hom^\bullet(F_{\uom},\mathbf{1})$ spanned by the $L_{\alpha,\eta}$ corresponding to these pairs. Then by Theorem \ref{thm_unitr}, 
$d(C_p)\subset C_{p+1}$, so
the $C_p$ are subcomplexes and they form a filtration of $\Hom^\bullet(F_{\uom},\mathbf{1})$. This allows one to define a spectral sequence as in \cite[Ch.\ 1, \S 4]{God}.
\begin{dfn}
    The \emph{path-dominance spectral sequence} for $\Hom^\bullet(F_{\uom},\mathbf{1})$, relative to $f$, is the spectral sequence associated to the above filtration.
\end{dfn}
The first page $E_1^{p,q}$ of this spectral sequence admits a general description. We will need a preliminary result. Let $\eta\in\{0,1\}^k$ and let $x_1,\dots, x_k$ be the corresponding Bruhat stroll. Let $\gamma_i:=x_i(\alpha_i)$ where $\alpha_i$ is the simple root corresponding to $s_i$. Then define $C_{\eta}$ to be the following complex:
    	    \begin{equation}
    	        \label{eq: full tensor decomposition}
		        \bigotimes_{\substack{\eta_i=U0 \\ e_i=+1}}(R\xrightarrow{id}\underline{R}) \otimes
		        \bigotimes_{\substack{\eta_i=U0 \\ e_i=-1}}(\underline{R}(-1)\xrightarrow{\gamma_i}R(1)) \otimes 
		        \bigotimes_{\substack{\eta_i=D0 \\ e_i=+1}}(R(-1)\xrightarrow{\gamma_i}\underline{R}(1)) \otimes
		        \bigotimes_{\substack{\eta_i=D0 \\ e_i=-1}}(\underline{R}\xrightarrow{id}R) 
    	    \end{equation}
In particular, when $\eta\in\mathcal{C}(\uom)$ there are no positive $U0$ or negative $D0$ so this is:
\begin{equation}
    \label{eq: critical tensor decomposition}
   \bigotimes_{\substack{\eta_i=U0 \\ e_i=-1}}(\underline{R}(-1)\xrightarrow{\gamma_i}R(1)) \otimes 	  \bigotimes_{\substack{\eta_i=D0 \\ e_i=+1}}(R(-1)\xrightarrow{\gamma_i}\underline{R}(1))    
\end{equation}
\begin{pro}\label{pro_classdiff}
	The row $E_0^{p,\bullet}$ of the path-dominance spectral sequence is homotopy equivalent to $\bigoplus_{\eta} C_\eta[p]$ where the sum is over $\eta\in \mathcal{C}(\uom)$ with $f(\eta)=p$.
\end{pro}
\begin{proof}
 Define any total order on 01-sequences which refines the path-dominance order. This defines a filtration on $\Hom^\bullet(F_{\uom},\mathbf{1})$ which refines the path-dominance filtration above. Each subquotient correspond to a 01-sequence $\eta$ and it is spanned by the dg light leaves $L_{\alpha,\eta}$ for all arrangements $\alpha$ that $\eta$ is compatible with. It easily follows from the proof of Theorem \ref{thm_unitr} that this subquotient is isomorphic to $C_\eta$. 
 
 This refined filtration induces a filtration of each $E_0^{p,\bullet}$ where the top subquotients correspond with critical 01-sequences $\eta$ with $f(\eta)=p$ and have no components of the differential between each other. All the below subquotients are contractible because there is at least a factor $R\overset{\id}{\to} R$ in their decomposition \eqref{eq: full tensor decomposition}. Hence the result follows.
\end{proof}
One can compute the cohomology of \eqref{eq: critical tensor decomposition} by means of the K\"unneth theorem and obtain the following.
Let $I_\eta$ be the ideal generated by the $\gamma_i$ in $R$ and let $d_\eta$ be the dimension of the span of the $\gamma_i$ in $\mathfrak{h}$. Let $c_\eta$ be the number of $0$'s in $\eta$ and $q_{\eta}=\sum_{\eta_i=0}e_i$. Then:
\begin{equation}
    \label{eq: cohomology of tensor product}
    H^{\bullet}(C_\eta) = \bigoplus_{i=d_\eta}^{c_\eta} R/I_\eta(c_\eta-2i)^{\oplus {c_\eta-d_\eta \choose i -d_\eta}}[-i+q_\eta]
\end{equation}
Then one can compute the first page of the path-dominance spectral sequence.
\begin{cor}\label{cor_firstpage}
    The rows of the first page of the path-dominance spectral sequence are:
    \begin{equation}
        \label{eq: description of first page}
        E_1^{p,\bullet}=\bigoplus_{f(\eta)=p} H^\bullet(C_\eta)[p] 
    \end{equation}
    the sum being over critical 01-sequences expressing the identity.
\end{cor}
\subsection{Cohomological range of morphisms}
We will use the results in the previous section to study the morphism space between standard Rouquier complexes. Let $x,y\in W$ and let $F_x$ and $F_y$ be the Rouquier complexes corresponding to the positive lifts of $x$ and $y$. 
\begin{thm}
    The morphism space $\Hom^\bullet(F_x,F_y)$ is zero unless $x\le y$ in the Bruhat order. In this case we have:
    \[
        \Hom^i(F_x,F_y)=0 \quad \text{unless $d\le i \le \ell(y)-\ell(x),$}
    \]
    where $d$ is the minimal dimension in $\mathfrak{h}$ of a span of roots $\alpha_{t_1},\dots \alpha_{t_k}$ associated to reflections $t_1,\dots t_k$ such that $x\to xt_1\to \dots \to xt_1\cdots t_k = y$ is a path in the Bruhat graph.
\end{thm}
\begin{proof}
    Consider reduced expressions $\underline{x}$ and $\underline{y}$ for $x$ and $y$ respectively. The morphism space in the statement is the same as the cohomology of $\Hom^\bullet (F_{\uom},\mathbf{1})$ where $\uom=\underline{x}\underline{y}^{-1}$. Notice that a critical 01-sequence $\eta$ for $\uom$ is forced to have all ones for the symbols corresponding to $\underline{x}$, because if there were zeros, the first would be decorated by $U$ by reducedness. If $\eta$ expresses the identity then there must be a subexpression of $\underline{y}$ which expresses $x$ and this is equivalent to $x\le y$. Then the first assertion follows from Corollary \ref{cor_firstpage}.

    From the above, one also gets that the maximal number of zeros in a critical 01-sequence $\eta$ is $c_\eta=\ell(y)-\ell(x)$ and $q_\eta=-c_\eta$ (all the zeros are negative). Hence the second statement follows from \eqref{eq: cohomology of tensor product} if we take the number of zeros as the path-dominant function.
\end{proof}
\begin{cor}
    If $W$ is a finite Weyl group and with standard realisation $\mathfrak{h}$, then 
    \[
        \Hom^i(F_x,F_y)=0 \quad \text{unless $d(x,y)\le i \le \ell(y)-\ell(x)$,}
    \]
    where $d(x,y)$ is the distance between $x$ and $y$ in the Bruhat graph. 
\end{cor}
\begin{proof}
    It suffices to notice that under the hypotheses, the roots associated to a minimal path in the Bruhat graph are linearly independent.
\end{proof}
\section{Extensions between Verma modules}
We will now give an application of the spectral sequence introduced in the previous section to category $\mathcal{O}$. As explained in the introduction, Rouquier complexes incarnate braid group actions on several categories, including $\mathcal{D}^{\mathrm{b}}(\mathcal{O}_0)$, the derived principal block of category $\mathcal{O}$. We will show that the above spectral sequence gives information about the extensions between Verma modules.
\subsection{Review of graded category $\mathcal{O}$}
Let $\mathfrak{g}$ be a complex semisimple Lie algebra and fix a triangular decomposition $\mathfrak{g}=\mathfrak{n}^-\oplus \mathfrak{h}\oplus\mathfrak{n}^+$. Let $\mathcal{O}$ be the Bernstein-Gelfand-Gelfand category of finitely generated $\mathfrak{h}$-semisimple $\mathfrak{n}^+$-locally nilpotent $\mathfrak{g}$-modules. Recall that to any $\lambda\in\mathfrak{h}^*$ one associates the Verma module:
\begin{equation}
    \Delta(\lambda):=\mathcal{U}(\mathfrak{g})\underset{\mathcal{U}(\mathfrak{b})}{\otimes} \mathbb{C}_\lambda,
\end{equation}
where $\mathcal{U}(\mathfrak{g})$ and $\mathcal{U}(\mathfrak{b})$ are the universal envelopping algebras associated, respectively,  to $\mathfrak{g}$ and to $\mathfrak{b}:=\mathfrak{h}\oplus\mathfrak{n}^+$, and $\mathbb{C}_\lambda$ is the one-dimensional $\mathfrak{b}$-module associated to $\lambda$. Simple modules of $\mathcal{O}$ are the heads of Verma modules and are denoted by $L(\lambda)$. 

Recall that the principal block $\mathcal{O}_0$ is the indecomposable direct summand of $\mathcal{O}$ generated by $L(0)$. This block contains a finite number of Verma modules and simple modules $\Delta_w:=\Delta(w\cdot 0)$ and $L_w:=L(w\cdot 0)$, where $w\in W$ is an element of the Weyl group $W$ associated to $\mathfrak{g}$ and $\mathfrak{h}$ and $w\cdot \lambda$ denotes the usual shifted action.
\subsection{Extensions between Verma modules}
It is an important open problem to determine the extensions between Verma modules. One reduces to the principal block $\mathcal{O}_0$ and the question becomes to determine the groups:
\begin{equation}
    \label{eq: ungraded extension}
    \mathrm{Ext}^i_{\mathcal{O}}(\Delta_x,\Delta_y).
\end{equation}
This problem has been studied extensively since the beginning of the theory. See, for example, the works of Delorme \cite{Delorme}, Gabber-Joseph \cite{GabberJoseph} and  Carlin \cite{Carlin1986}. It was expected that the dimension of \eqref{eq: graded extensions} should be prescribed by the $R$-polynomial $R_{x,y}$:
\begin{equation}
    \label{eq:expected_dimensions}
    \dim \mathrm{Ext}^k_{\mathcal{O}}(\Delta_x,\Delta_y )=
        |R_{x,y}^{(k)}|
\end{equation}
However Boe \cite{Boe} proved that ``additional'' extensions exist.

One can refine the question by passing to the graded version $\mathcal{O}_0^{\mathbb{Z}}$ of $\mathcal{O}_0$ and investigate the groups:
\begin{equation}
    \label{eq: graded extensions}
    \mathrm{Ext}^k_{\mathcal{O}_0^{\mathbb{Z}}}(\Delta_x\langle i\rangle, \Delta_y).
\end{equation}
In a recent paper \cite{MazKo}, Ko and Mazorchuk observed that the information provided by the coefficients of $R_{x,y}$ is actually the alternate sum:
\begin{equation}
    \label{eq: information provided by R polynomials}
    R_{x,y}^{(i)}=\sum_{k\ge 0}(-1)^k\dim \mathrm{Ext}_{\mathcal{O}_0^{\mathbb{Z}}}^k( \Delta_x\langle i \rangle,\Delta_y).
\end{equation}
They also described bounds for the dimension of the extensions \eqref{eq: graded extensions} and described the region of pairs $(k,i)$ such that extensions can exist. This is the triangle defined by: $0\le k\le \ell(x)-\ell(y)$, $k\le i$ and $k\ge 2i-\ell(x)+\ell(y)$, with the top and left side removed, except the top and bottom vertex (an example is depicted in the introduction: the white dots are removed). Additional extension are precisely those corresponding to points in the interior of the triangle.

We will give an interpretation of \eqref{eq: information provided by R polynomials} and of this region in terms of our spectral sequence and we will deduce a criterion to detect additional extension.
\subsection{Soergel's approach to category $\mathcal{O}$}
Soergel \cite{Soe_HC} observed that one can realise the projective objects of $\mathcal{O}_0$ as modules over the coinvariant algebra $R/R_+^W$. Furthermore, after embedding the latter into just modules over $R$, one can characterise the modules that correspond to projectives in $\mathcal{O}_0$. These are called \textit{Soergel modules}. In terms of the general version of the Hecke category introduced in \S \ref{sec: Hecke category}:
\begin{dfn}
    The category $\overline{\mathscr{H}}$ has the same objects as $\mathscr{H}$ but morphism spaces are:
    \begin{equation}
        \label{eq: morphisms in Soergel modules}
        \Hom_{\overline{\mathscr{H}}}(\overline{B},\overline{B}'):=\Bbbk\otimes_R \Hom_{\mathscr{H}}(B,B')
    \end{equation}
    via the map $R\to \Bbbk$ that sends $\mathfrak{h}^*$ to 0.
\end{dfn}
For Bott-Samelson modules (i.e., images of Bott-Samelson bimodules), the diagrammatic computations work in the same way except any morphism with a homogeneous polynomial of positive degree on the left is zero.

The category $\overline{\mathscr{H}}$ is still graded, but now morphism spaces are only $R$-modules and not bimodules. This is a module category over the Hecke category via $\overline{B}\cdot B':=\overline{BB'}$ and the obvious maps on morphisms. Theorem \ref{thm: categorification} gives immediately the following.
\begin{thm}
    \label{thm: categorification module}
    Let $\Bbbk$ be a complete local ring. The map $K_{\oplus}(\overline{\mathscr{H}})\to \mathbf{H}$ uniquely determined by $[\overline{B}_s]\mapsto \delta_s+v$ is an isomorphism. Furthermore the morphism spaces $\Hom_{\overline{\mathscr{H}}}(\overline{B}_1,\overline{B}_2)$ are free as graded $\Bbbk$-modules and the graded dimension is given by:
    \begin{equation}
        \label{eq: graded dimension over k}
        \mathrm{grdim}_\Bbbk \Hom(\overline{B}_1,\overline{B}_2) = \beta([\overline{B}_1],[\overline{B}_2]).
    \end{equation}
\end{thm}
\begin{rmk}
    The category $\overline{\mathscr{H}}$ is then a categorification of the regular representation of the Hecke algebra.
\end{rmk}
 Let $\overline{\mathscr{H}}^{\circ}$ be the corresponding ungraded version. Then Soergel's results imply:
\begin{equation}
    \label{eq: Soergel description of derived O}
        \mathcal{K}^{\mathrm{b}}(\overline{\mathscr{H}}^\circ)\cong \mathcal{K}^{\mathrm{b}}(\mathrm{Proj-}\mathcal{O}_0)\cong \mathcal{D}^{\mathrm{b}}(\mathcal{O}_0)
\end{equation}
Hence the Soergel theoretic version of $\DbOZ$ is just $\Kb(\overline{\mathscr{H)}}$.
Then the braid group action on $\mathcal{D}^{\mathrm{b}}(\mathcal{O}_0)$ comes from the natural action of $\Kb(\mathscr{H})$ on $\Kb(\overline{\mathscr{H}})$ by restricting to the monoidal subcategory of Rouquier complexes.

Furthermore, the standard $t$-structure on $\DbOZ$ corresponds to the $t$-structure defined by Achar and Riche \cite{AchRic_modII} on $\Kb(\overline{\mathscr{H}})$ and $\langle 1 \rangle$ corresponds to the Tate twist $\langle 1\rangle:=[-1](1)$. The heart of the $t$-structure contains the images $\overline{F}_x$ of all the standard Rouquier complexes and under the equivalence with $\mathcal{O}_0^{\mathbb{Z}}$ we have:
\begin{equation}
    \label{eq: standard and costandard correspondence}
    \Delta_x \quad \leftrightarrow \quad \overline{F}_{w_0x^{-1}}
\end{equation}
\subsection{Adaptation of the spectral sequence}
The dg description of Rouquier complexes carries over to their images in $\Kb(\overline{\mathscr{H}})$.
For Bott-Samelson modules the light leaves produce a a basis of for morphism spaces $\Hom(\overline{B}_{\underline{x}},\mathbf{1})$, this time over $\Bbbk$. So the dg light leaves produce a basis over $\Bbbk$ for the dg morphism spaces $\Hom^\bullet(\overline{F}_{\uom},\mathbf{1})$. 
Given a path-dominant function $f$, one can define a spectral sequence $\overline{E}^{p,q}_r$ in the same way as in \ref{sec: spectral sequence} and one can deduce its properties from those of the original one.
More precisely, if $\eta\in \mathcal{C}(\uom)$ then let:
\begin{equation}
    \label{eq: module version of tensor decomposition}
    \overline{C}_\eta:=\bigotimes_{\eta_i=U0}(\underline{\Bbbk}(-1)\overset{0}{\to}\Bbbk(1))\otimes 
    \bigotimes_{\eta_i=D0}(\Bbbk(-1)\overset{0}{\to}\underline{\Bbbk}(1))
\end{equation}
\begin{cor}
    The row $\overline{E}^{p,\bullet}_0$ of the zero page is homotopy equivalent to $\bigoplus_{\eta} \overline{C}_{\eta}[p]$ where the sum is over $\eta\in \mathcal{C}(\uom)$ with $f(\eta)=p$.
\end{cor}
\begin{proof}
    This follows in the same way as Proposition \ref{pro_classdiff} by noticing that the $\gamma_i$ become $0$ in the module category.
\end{proof}
Because all the differentials in $\overline{C}_\eta$ are zero we get:
\begin{cor}
    \label{cor: rows of first page for soergel modules}
    The rows of the first page are $\overline{E}^{p,\bullet}_1\cong \bigoplus_\eta \overline{C}_\eta[p]$ where the sum is over $\eta\in \mathcal{C}(\uom)$ with $f(\eta)=p$.
\end{cor}
\subsection{Extension groups}
The extension groups \eqref{eq: graded extensions} can be obtained as:
\begin{equation}
    \label{eq: extensions as Hom in derived category}
    \mathrm{Ext}^k_{\mathcal{O}_0^{\mathbb{Z}}}(\Delta_x\langle i \rangle,\Delta_y)=\Hom_{\DbOZ}(\Delta_x\langle i \rangle, \Delta_y [k]).
\end{equation}
 Hence one can compute these extensions as:
\begin{equation}
    \label{eq: extensions as Hom in homotopy category}
    \mathrm{Ext}^k_{\mathcal{O}_0^{\mathbb{Z}}}(\Delta_x\langle i \rangle,\Delta_y)=\Hom_{\Kb(\overline{\mathscr{H}})}(\overline{F}_{w_0 x^{-1}}\langle i \rangle, \overline{F}_{w_0y^{-1}} [k]).
\end{equation}
By adjunction, one can then compute:
\begin{equation}
    \Hom_{\Kb(\overline{\mathscr{H}})}(\overline{F}_{yx^{-1}}\langle i \rangle,\mathbf{1}[k])
\end{equation}
Let $\underline{x}$ and $\underline{y}$ be reduced expressions for $x$ and $y$ respectively. 
\begin{lem}
    In any critical 01-sequence for $\uom=\underline{y}\underline{x}^{-1}$ the symbols for the letters in $\underline{y}$ are all 1.
\end{lem}
\begin{proof}
    Suppose one of them is 0. Because $\underline{y}$ is reduced, the first such 0 would be decorated with $U$, so the 01-sequence would not be critical.
\end{proof}
Critical 01-sequences of $\uom$ correspond then to distinguished 01-sequences for $\underline{y}$ expressing $x$. One can define the \textit{standard} path-dominant function as in Example \ref{ex: path-dominant functions} (ii). More precisely, if $h\in \{0,1\}$ such that $\ell(x)-\ell(y)\equiv h\mod 2$, then for every $\eta$ with $g(\eta)$ 0s, let $f(\eta)=(g(\eta)-h)/2$. Hence we can define the spectral sequence $\overline{E}^{p,q}_r$ corresponding with $f$.
\begin{cor}
    \label{cor: first page for extensions}
    The rows of the first page of the spectral sequence above are:
    \begin{equation}
    \label{eq: rows of first page for extension}
        \overline{E}^{p,\bullet}_1\cong \bigoplus_{f(\eta)=p}(\Bbbk(-1)\oplus \Bbbk[-1](1))^{\otimes (2p+h)} [p]
    \end{equation}
    So the Poincar\'e polynomial of $\overline{E}^{\bullet,\bullet}_1$ is
        $P(s,t)=s^{-\frac{h}{2}}t^{\frac{h}{2}}\widetilde{R}_{x,y}(s^{\frac{1}{2}}(t^{-\frac{1}{2}}+t^{\frac{1}{2}}))$.
\end{cor}
\begin{proof}
    The first statement is an immediate consequence of Corollary \ref{cor: rows of first page for soergel modules}.

    Then notice that each summand of \eqref{eq: rows of first page for extension} contributes to $P(s,t)$ by:
    \[
    	s^pt^{-p}(1+t)^{2p+h}=s^{\frac{g(\eta)-h}{2}}t^{-\frac{g(\eta)-h}{2}}(1+t)^{g(\eta)}=
    	s^{-\frac{h}{2}}t^{\frac{h}{2}}(s^{\frac{1}{2}}(t^{-\frac{1}{2}}+t^{\frac{1}{2}}))^{g(\eta)}.
    \]
    Hence the second statement follows from Proposition \ref{pro: classical properties of R - polynomials} (ii).
\end{proof}
\begin{pro}
    The spectral sequence $\overline{E}^{p,q}_r$ converges at the second page and $\overline{E}^{p,q}_2\cong \mathrm{Ext}^k_{\mathcal{O}^{\mathbb{Z}}_0}(\Delta_x\langle i \rangle,\Delta_y)$ if $k=p-q+h$ and $i=h-2q$.
\end{pro}
\begin{proof}
For the convergence statements it is sufficient to observe that the differential of the second and the following pages must be 0 for degree reasons.

Then notice that $\overline{E}^{p,q}_2$ is a sum of copies of $\Bbbk(2q-h)$, hence it contributes to $\Hom(F_{yx^{-1}},\mathbf{1}(h-2q)[p+q])=\mathrm{Ext}^{p-q+h}(\Delta_x\langle h-2q\rangle,\Delta_y)$. Hence the second statement follows.
\end{proof}
Arranging the possible contributions on the $(k,i)$ plane then gives the same triangle as in \cite{MazKo}. Indeed for the spectral sequence we have that $\overline{E}_1^{p,q}=0$ unless $0\le p \le (\ell(x)-\ell(y)-h)/2$ and $-p\le q\le p+h$, which is equivalent to $0\le k \le \ell(x)-\ell(y)$ and $2k-\ell(x)+\ell(y)\le i\le k$.
\subsection{A criterion for additional extensions}
The expected extensions correspond to the $\overline{E}_2^{p,q}$ with $p$ maximal, that is $p=(\ell(x)-\ell(y)-h)/2$. Additional extension correspond to any other value of $p$.

Now, the differential of the first page is of the form $\overline{E}^{p,q}_1\to \overline{E}_1^{p+1,q}$. So fix $q$ and consider the increasing sums:
\begin{equation}
    \label{eq: alternating sums on the spectral sequence}
    S_r^q:=\sum_{p=0}^r (-1)^p\dim \overline{E}_1^{p,q}.
\end{equation}
Then if the sequence $(S_0^q,S_1^q,\dots)$ does not alternate in sign, this means that some $\overline{E}^{p,q}_2\neq 0$ for non maximal values of $p$, that is, there are additional extensions for that value of $q$ (i.e., for $i=h-2q$).
\bibliographystyle{amsplain}
\bibliography{bibliografia}

\end{document}